\documentclass[a4,10pt]{article}
\usepackage{amsmath,amscd,amssymb}
\usepackage[margin=.8in]{geometry}

\newcommand{\nbiga}{\mathcal{A}}

\newcommand{\nbigc}{\mathcal{C}}
\newcommand{\nbigd}{\mathcal{D}}

\newcommand{\nbigf}{\mathcal{F}}
\newcommand{\nbigg}{\mathcal{G}}

\newcommand{\nbigl}{\mathcal{L}}
\newcommand{\nbigm}{\mathcal{M}}

\newcommand{\nbigo}{\mathcal{O}}

\newcommand{\nbigr}{\mathcal{R}}

\newcommand{\nbigt}{\mathcal{T}}

\newcommand{\nbigx}{\mathcal{X}}
\newcommand{\nbigy}{\mathcal{Y}}

\newcommand{\proj}{\mathbb{P}}
\newcommand{\seisuu}{{\mathbb Z}}

\newcommand{\cnum}{{\mathbb C}}
\newcommand{\real}{{\mathbb R}}

\newcommand{\gbigc}{\mathfrak C}
\newcommand{\gbigd}{\mathfrak D}

\newcommand{\gminib}{\mathfrak b}

\newcommand{\vecS}{{\boldsymbol S}}

\newcommand{\larr}{\leftarrow}

\newcommand{\lrarr}{\longrightarrow}

\newcommand{\pf}{{\bf Proof}\hspace{.1in}}
\newcommand{\qed}{\mbox{\rule{1.2mm}{3mm}}}

\def\Image{\mathop{\rm Im}\nolimits}

\def\Gr{\mathop{\rm Gr}\nolimits}

\def\Tot{\mathop{\rm Tot}\nolimits}

\def\rank{\mathop{\rm rank}\nolimits}

\def\Ker{\mathop{\rm Ker}\nolimits}

\def\Gr{\mathop{\rm Gr}\nolimits}
\def\Sym{\mathop{\rm Sym}\nolimits}

\def\id{\mathop{\rm id}\nolimits}

\newcommand{\del}{\partial}
\newcommand{\delbar}{\overline{\del}}

\newcommand{\nhom}{{\mathcal Hom}}

\newcommand{\barlambda}{\overline{\lambda}}
\newcommand{\lambdabar}{\barlambda}

\def\reg{\mathop{\rm reg}\nolimits}

\newcommand{\deldel}{\eth}

\newcommand{\distribution}{\gbigd\gminib}

\newcommand{\Rtilde}{\widetilde{R}}

\newcommand{\nbigrtilde}{\widetilde{\nbigr}}

\newcommand{\nbigmtilde}{\widetilde{\nbigm}}

\def\MT{\mathop{\rm MT}\nolimits}

\def\MTM{\mathop{\rm MTM}\nolimits}

\def\DR{\mathop{\rm DR}\nolimits}

\def\Ch{\mathop{\rm Ch}\nolimits}

\def\integral{\mathop{\rm int}\nolimits}
\def\pt{\mathop{\rm pt}\nolimits}

\def\alg{\mathop{\rm alg}\nolimits}

\def\Cr{\mathop{\rm Cr}\nolimits}

\newcommand{\Omegatilde}{\widetilde{\Omega}}

\newcommand{\iotatilde}{\widetilde{\iota}}

\newcommand{\alphatilde}{\widetilde{\alpha}}

\newcommand{\otimeshat}{\widehat{\otimes}}

\newtheorem{thm}{Theorem}[section]
\newtheorem{cor}[thm]{Corollary}

\newtheorem{rem}[thm]{Remark}
\newtheorem{lem}[thm]{Lemma}
\newtheorem{prop}[thm]{Proposition}

\newtheorem{question}[thm]{Question}

\begin{document}

\title{A generalization of Barannikov-Kontsevich theorem}

\author{Takuro Mochizuki\thanks{Research Institute for Mathematical Sciences, Kyoto University, Kyoto 606-8512, Japan, takuro@kurims.kyoto-u.ac.jp}}
\date{}
\maketitle

\begin{abstract}
We study the twisted de Rham complex associated with
a holomorphic function on a K\"ahler manifold
whose critical point set is compact.
We prove
the $E_1$-degeneration of
the Hodge-to-de Rham spectral sequence.
It is a generalization of Barannikov-Kontsevich Theorem. 

\vspace{.1in}
\noindent
MSC: 
\end{abstract}

\section{Introduction}

\subsection{Twisted de Rham complexes}

Let $X$ be a complex manifold.
Let $f$ be a holomorphic function on $X$
such that the set of critical points $\Cr(f)$ is compact.

Let $\Omega _X^k$ denote the sheaf of
holomorphic $k$-forms on $X$.
There exists the exterior derivative
$d:\Omega^{k}_X\to \Omega^{k+1}_X$.
The exterior product of $df$ induces morphisms
$df:\Omega^k_X\to \Omega^{k+1}_X$,
i.e.,
$\tau\longmapsto df\wedge\tau$.

Let $\lambda$ be a variable.
Let $\Omega_X^k[\![\lambda]\!]$
denote the sheaf of formal power series
with $\Omega_X^k$-coefficients.
Namely,
for any open subset $U$ of $X$,
let $\Omega_X^k[\![\lambda]\!](U)$
denote the space of 
formal power series
$\sum_{j\geq 0} \tau_j\lambda^j$,
where $\tau_j\in\Omega_X^k(U)$.
We obtain the differential
$\lambda d+df:
\Omega_X^{k}[\![\lambda]\!]
\to
\Omega_X^{k+1}[\![\lambda]\!]$
determined by
$\tau\mapsto \lambda d\tau+df\wedge \tau$.
It satisfies $(\lambda d+df)^2=0$.
The complex
$\Omega_X^{\bullet}[\![\lambda]\!]_f=
(\Omega_X^{\bullet}[\![\lambda]\!],\lambda d+df)$
is called the twisted de Rham complex,
or the formal twisted de Rham complex
when we emphasize to consider the formal series.
We obtain the cohomology group
$H^{\ast}(X,\Omega^{\bullet}_X[\![\lambda]\!]_f)$.
It is naturally a $\cnum[\![\lambda]\!]$-module,
where $\cnum[\![\lambda]\!]$ denotes the ring of formal power
series with $\cnum$-coefficients.

There are several different versions of
twisted de Rham complex.
For example,
if $X$ and $f$ are algebraic,
let $\Omega^{k}_{X^{\alg}}$ denote the sheaf of
algebraic $k$-forms
on the algebraic variety $X^{\alg}$
with Zariski topology,
and 
it is also natural to consider
the sheaf 
$\Omega_{X^{\alg}}^k[\lambda]$
of polynomials
with $\Omega_{X^{\alg}}^k$-coefficients
on $X^{\alg}$.
We obtain the algebraic twisted de Rham complex
$\Omega^{\bullet}_{X^{\alg}}[\lambda]_f=
(\Omega^{\bullet}_{X^{\alg}}[\lambda],\lambda d+df)$
and the cohomology group
$H^{\ast}(X^{\alg},\Omega^{\bullet}_{X^{\alg}}[\lambda]_f)$.

\subsection{A basic question}

By setting 
$\nbigf_{\lambda}^j(\Omega^{\bullet}_X[\![\lambda]\!]_f)
=\lambda^j\Omega^{\bullet}_X[\![\lambda]\!]_f$
for any non-negative integer $j$,
we obtain the filtered complex
$\nbigf_{\lambda}\Omega^{\bullet}_X[\![\lambda]\!]_f$.
The associated complexes
$\Gr^j_{\nbigf_{\lambda}}(\Omega^{\bullet}_X[\![\lambda]\!]_f)
=\nbigf_{\lambda}^j/\nbigf_{\lambda}^{j+1}$ $(j\geq 0)$
are isomorphic to
the complex $(\Omega^{\bullet}_X,df)$.
Note that the cohomological support of
$(\Omega^{\bullet}_X,df)$ is contained in $\Cr(f)$.
Because $\Cr(f)$ is compact,
$H^{\ast}(X,(\Omega^{\bullet}_X,df))$ is
a finite dimensional complex vector space.
There is the spectral sequence
associated with the filtered complex,
for which 
\[
E_1^{p,q}=
H^{p+q}\bigl(X,
\Gr^{p}_{\nbigf_{\lambda}}
(\Omega^{\bullet}_X[\![\lambda]\!]_f)
\bigr)
=H^{p+q}(X,(\Omega_X^{\bullet},df)).
\]
We study the following question.

\begin{question}
\label{question;25.6.26.1}
Is the spectral sequence
for $\nbigf_{\lambda}\Omega^{\bullet}_X[\![\lambda]\!]_f$
degenerates
at the $E_1$-level?
\end{question}

We can rephrase the condition in several ways.
\begin{lem}
\label{lem;25.6.29.20}
The spectral sequence degenerates at the $E_1$-level
if and only if the following equivalent conditions are satisfied.
\begin{itemize}
 \item For any $\ell$ and $j$,
      the exact sequences
       $0\lrarr
\Gr_{\nbigf_{\lambda}}^j\lrarr
\Omega_X^{\bullet}[\![\lambda]\!]_f/\nbigf_{\lambda}^{j+1}\lrarr
       \Omega_X^{\bullet}[\![\lambda]\!]_f/\nbigf_{\lambda}^{j}\lrarr 0$
       induce the following exact sequences:
\[
 0\lrarr
 H^{\ell}\bigl(
 X,\Gr_{\nbigf_{\lambda}}^j\bigl(\Omega_X^{\bullet}[\![\lambda]\!]_f\bigr)
 \bigr)
 \lrarr
  H^{\ell}\bigl(
 X,\Omega_X^{\bullet}[\![\lambda]\!]_f/\nbigf_{\lambda}^{j+1}
 \bigr)
 \lrarr \\
 H^{\ell}\bigl(
 X,\Omega^{\bullet}_X[\![\lambda]\!]_f/\nbigf_{\lambda}^{j}
 \bigr)
 \lrarr 0.
\]      
 \item 
 For any $\ell$,
  $H^{\ell}(X,\Omega^{\bullet}_X[\![\lambda]\!]_f)$
      are isomorphic to
      $H^{\ell}(X,(\Omega^{\bullet}_X,df))\otimes\cnum[\![\lambda]\!]$
      as $\cnum[\![\lambda]\!]$-modules.
In particular,
the formal twisted de Rham cohomology groups
$H^{\ell}(X,\Omega^{\bullet}_X[\![\lambda]\!]_f)$
are free
$\cnum[\![\lambda]\!]$-modules
       of finite rank.
\end{itemize}
\end{lem}

\begin{rem}
It is easy to see that the second condition implies
the first.
For the convenience of readers,
we shall recall that the $E_1$-degeneration
is equivalent to the first condition
in {\rm\S\ref{subsection;25.8.25.10}}.
We shall later explain that
the first condition implies the second condition
in a generalized context. 
(See Corollary {\rm\ref{cor;25.7.2.1}}.)
\end{rem}

\subsubsection{Classical results related to the question}

There are two classical results
related to this question.
One appeared in the classical Hodge theory.
As Deligne observed,
if $X$ is projective and $f=0$,
the desired $E_1$-degeneration
follows from the Hodge decomposition of $H^{\ast}(X,\cnum)$.
Deligne generalized it to
more general algebraic varieties by using mixed Hodge theory.
This $E_1$-degeneration and its generalizations
are not only deep results in the Hodge theory,
but also useful for various applications
including some vanishing theorems.
Deligne and Illusie gave
an alternative proof of the $E_1$-degeneration
using the reduction to the positive characteristic.

The other appeared in the singularity theory,
in particular, the study of Brieskorn lattices.
If $f$ has
only one critical point,
the $E_1$-degeneration holds
because
the $\ell$-th cohomology group
of $(\Omega^{\bullet}_X,df)$
is $0$
unless $\ell$ equals $\dim X$.
In this case,
the $E_1$-degeneration is an important starting point of
the deep theory of primitive forms of Kyoji Saito.
(See \cite{Sabbah-Brieskorn}
and \cite{Saito-Takahashi} for more backgrounds.)

\subsection{Barannikov-Kontsevich Theorems and variations}
\label{subsection;25.6.26.2}

The modern study of twisted de Rham complexes
was opened by
the celebrated theorem
of Barannikov and Kontsevich.
\begin{thm}[Barannikov-Kontsevich]
\label{thm;25.6.29.31}
The $E_1$-degeneration
for $\nbigf_{\lambda}\Omega^{\bullet}_X[\![\lambda]\!]_f$
holds 
in the case where 
$X$ and $f$ are quasi-projective.
\end{thm}

This is a fundamental theorem
in the study
of the holomorphic Landau-Ginzburg model
of the mirror symmetry.
For example,
it is essential
in the proof of
smoothness of
some moduli spaces
associated with
Landau-Ginzburg models.
(See \cite{Katzarkov-Kontsevich-Pantev-2017}.)

The original proof of
Barannikov and Kontsevich
was given by the harmonic analysis
on the basis of a generalization of the K\"ahler identity to this context.
\footnote{
In the first version of this paper,
due to my mistaken assumption,
I wrote that 
``The original proof of Barannikov and Kontsevich
was given by a generalization of
the method of Deligne and Illusie.''
I thank Maxim Kontsevich for pointing out it.}
Indeed, 
Barannikov and Kontsevich proved the following
theorem for the algebraic version of
the twisted de Rham complexes,
which implies Theorem \ref{thm;25.6.29.31}.

\begin{thm}[Barannikov-Kontsevich]
\label{thm;25.6.29.30}
Suppose that $f:X\to \cnum$ is a projective morphism of
algebraic varieties.
Then, 
\[
  \dim H^{j}(X^{\alg},(\Omega^{\bullet}_{X^{\alg}},df))
=\dim H^{j}(X^{\alg},(\Omega^{\bullet}_{X^{\alg}},d+df))
\]
holds for any $j$ and for any complex number $\lambda$.
\hfill\qed
\end{thm}

Theorem \ref{thm;25.6.29.30} implies that
$H^j(X^{\alg},\Omega_{X^{\alg}}^{\bullet}[\lambda]_f)$
are free $\cnum[\lambda]$-modules,
and that 
the $E_1$-degeneration
of the spectral sequence
for the filtered complex
$\nbigf_{\lambda}\Omega^{\bullet}_{X^{\alg}}[\lambda]_f$.
In the setting of Theorem \ref{thm;25.6.29.31},
there exists a projective morphism
of algebraic varieties
$F:Y\to \cnum$ 
with an open embedding
$\iota:X\to Y$ such that $f=F\circ\iota$.
Under the assumption that $\Cr(f)$ is compact,
the set $\Cr(F)$ of critical points of $F$
is decomposed as
$\Cr(F)=\Cr(f)\sqcup\bigl(
\Cr(F)\cap(Y\setminus X)\bigr)$.
Then, we obtain the $E_1$-degeneration for
$\nbigf_{\lambda}\Omega^{\bullet}_X[\![\lambda]\!]_f$
from the $E_1$-degeneration for
$\nbigf_{\lambda}\Omega^{\bullet}_{Y^{\alg}}[\lambda]_f$.

The theorem of Barannikov-Kontsevich
for the algebraic twisted de Rham complexes
(Theorem \ref{thm;25.6.29.30})
has attracted many mathematicians
because of its significance
in the non-commutative Hodge theory
(see \cite{Katzarkov-Kontsevich-Pantev-2017}),
and because the theorem and its generalization are
deeply related with various fields of mathematics.
Indeed,
alternative proofs
for Theorem \ref{thm;25.6.29.30}
with different methods have been found
by Sabbah \cite{Sabbah-twisted-de-Rham}
using Hodge modules and microlocalization,
and by Ogus and Vologodsky \cite{Ogus-Vologodsky}
using their non-abelian Hodge correspondence
in positive characteristic.
Later,
Arinkin, C\u{a}ld\u{a}raru and Hablicsek
\cite{Arinkin-Caldararu-Hablicsek}
revisited it in their study of
Deligne-Illusie method from
the viewpoint of derived algebraic geometry.
Sabbah also studied generalizations
to the case where $f$ is not necessarily projective
but cohomologically tame \cite{Sabbah-tame-polynomial}.
See \cite{Esnault-Sabbah-Yu} and
\cite{Mochizuki-Kontsevich-complexes}
for a generalization to the Kontsevich complexes.

The $L^2$-analogue 
of the twisted de Rham complex associated with $(X,f)$
has been studied
also in \cite{Fan,Li-Wen}.
In particular,
Li and Wen \cite{Li-Wen} studied the case
where $X$ has a complete K\"ahler metric
with bounded curvature,
and $f$ is strongly elliptic,
which is
a kind of non-degeneracy condition
at infinity.
They established
an analogue of Theorem \ref{thm;25.6.29.30}
in this context.

\subsection{Main result}

In this paper,
we shall study a generalization of
Theorem \ref{thm;25.6.29.31}.
It is an affirmative answer to a question
asked by Kontsevich to the author.

\begin{thm}
\label{thm;23.3.29.1}
If $X$ is K\"ahler,
the $E_1$-degeneration 
for $\nbigf_{\lambda}\Omega_X^{\bullet}[\![\lambda]\!]_f$
holds.
\end{thm}

Note that $X$ can be a small neighbourhood of $\Cr(f)$,
and that we do not need any assumption
on the behaviour of $(X,f)$ at infinity.
It is our purpose to show that
a local assumption around $\Cr(f)$ is enough for
the $E_1$-degeneration for
the formal twisted de Rham complex,
though global assumptions are useful
to obtain stronger consequences
as in Theorem \ref{thm;25.6.29.30}.

\paragraph{Acknowledgement}

I would like to thank Maxim Kontsevich for
what he has brought to mathematics.
This study is motivated by one of his questions.
I also appreciate him for pointing out a historical inaccuracy
in a previous version.
Theorem \ref{thm;23.3.29.1}
is one of three topics
in  the author's talk in the conference
``Mathematics on the Crossroad of Centuries''
held in 2024 September.
The others will be explained elsewhere.
A preliminary version of this paper was prepared for
the conference in Rikkyo University in 2023 January.
I heartily thank the organizers for the opportunity of the talks.

I thank Claude Sabbah for helpful discussions
and for his kindness.

I am partially supported by
the Grant-in-Aid for Scientific Research (A) (No. 21H04429),
the Grant-in-Aid for Scientific Research (A) (No. 22H00094),
the Grant-in-Aid for Scientific Research (A) (No. 23H00083),
and the Grant-in-Aid for Scientific Research (C) (No. 20K03609),
Japan Society for the Promotion of Science.
I am also partially supported by the Research Institute for Mathematical
Sciences, an International Joint Usage/Research Center located in Kyoto
University.

\section{Integrable mixed twistor $\nbigd$-modules}

\subsection{Mixed twistor $\nbigd$-modules}

The theory of twistor $\nbigd$-modules
has been developed
in \cite{Sabbah-pure-twistor,Sabbah-wild-twistor}
and \cite{mochi2,Mochizuki-wild,Mochizuki-MTM}
as a twistor version of Hodge modules
\cite{saito1,saito2}
inspired by the principle called Simpson's Meta theorem \cite{s3}.

\subsubsection{$\nbigr_X$-triples}

For any complex manifold $X$,
we set $\nbigx=\cnum_{\lambda}\times X$.
Let $p_{\lambda}:\nbigx\to X$ denote the projection.
Let $\nbigd_{\nbigx}$ denote the sheaf of
holomorphic linear differential operators on $\nbigx$,
and let $\Theta_X$ denote the tangent sheaf of $X$.
We obtain the sheaf of subalgebras 
$\nbigr_X\subset\nbigd_{\nbigx}$
generated by
$\lambda \cdot (p_{\lambda}^{\ast}\Theta_X)$ over $\nbigo_{\nbigx}$.
If $X$ is an open subset in $\cnum^n$,
we have
$\nbigr_X=
\nbigo_{\nbigx}\langle\lambda\del_1,\ldots,\lambda\del_n\rangle$.

We set
$\vecS=\bigl\{\lambda\in\cnum\,\big|\,|\lambda|=1\bigr\}$.
Let $\sigma:\vecS\to\vecS$ be defined by
$\sigma(\lambda)=-\lambda$.
Let $\distribution_{\vecS\times X/\vecS}$
denote the sheaf of distributions
on $\vecS\times X$ which are continuous with respect to $\vecS$.
(See \cite[\S0.5]{Sabbah-pure-twistor}.)
It is naturally
an $\nbigr_{X|\vecS\times X}\otimes_{\nbigo_{\cnum_{\lambda}|\vecS}}
\sigma^{-1}(\nbigr_{X|\vecS\times X})$-module
by the action
$(P_1\otimes\sigma^{-1}(P_2))\cdot \tau
=P_1\overline{\sigma^{-1}(P_2)}\tau$.
A sesqui-linear pairing of $\nbigr_X$-modules $\nbigm'$ and $\nbigm''$
is a morphism of
$\nbigr_{X|\vecS\times X}\otimes_{\nbigo_{\cnum_{\lambda}|\vecS}}
\sigma^{-1}(\nbigr_{X|\vecS\times X})$-modules
$\nbigm'_{|\vecS\times X}\otimes_{\nbigo_{\cnum_{\lambda}|\vecS}}
\sigma^{-1}(\nbigm''_{|\vecS\times X})
\to
\distribution_{\vecS\times X/\vecS}$.
Such a tuple $(\nbigm',\nbigm'',C)$
is called an $\nbigr_X$-triple.

For $\nbigr_X$-triples
$\nbigt_i=(\nbigm_i',\nbigm_i'',C_i)$ $(i=1,2)$,
a morphism 
$\nbigt_1\to\nbigt_2$
is defined to be a pair of $\nbigr_X$-homomorphisms
$\varphi':\nbigm_2'\to\nbigm_1'$
and $\varphi'':\nbigm_1''\to\nbigm_2''$
such that
$C_1\circ (\varphi'\times\id)
=C_2\circ(\id\times\varphi'')$.
The category of $\nbigr_X$-triples
is an abelian category.

For an increasing filtration $W$ of $\nbigt$ 
in the category of $\nbigr_X$-triples,
we have the increasing filtrations
$W(\nbigm')$ and $W(\nbigm'')$
such that
$W_j(\nbigt)=\bigl(
 \nbigm'/W_{-j-1}(\nbigm'),
 W_j(\nbigm''),C_j
\bigr)$,
where $C_j$ denote the induced sesqui-linear pairings.

\subsubsection{Direct image of $\nbigr_X$-triples}

Let $F:X\lrarr Y$ be a morphism of complex manifolds.
We set $\omega_{\nbigx}:=\lambda^{-\dim X}p_X^{\ast}\omega_{X}$,
where $\omega_X$ denotes the canonical line bundle of $X$.
\index{sheaf $\omega_{\nbigx}$}
\index{sheaf $\omega_X$}
Similarly, we set
$\omega_{\nbigy}:=\lambda^{-\dim Y}p_Y^{\ast}\omega_Y$.
We set
$\nbigr_{Y\larr X}:=
 \omega_{\nbigx}\otimes_{F^{-1}(\nbigo_{\nbigy})}
 F^{-1}(\nbigr_{Y}\otimes\omega_{\nbigy}^{-1})$. 
\index{sheaf $\nbigr_{Y\larr X}$}
For any $\nbigr_{X}$-module $\nbigm$,
we obtain the following $\nbigr_{Y}$-modules:
\[
 F^j_{\dagger}(\nbigm):=
 R^j(\id_{\cnum}\times F)_{!}
 \bigl(
  \nbigr_{Y\larr X}\otimes^L_{\nbigr_{X}}
  \nbigm
  \bigr)
  \quad
  (j\in\seisuu).
\]

For any $\nbigr_X$-triple $\nbigt=(\nbigm',\nbigm'',C)$,
we obtain the $\nbigr_X$-triples
$F_{\dagger}^j(\nbigt)
=(F_{\dagger}^{-j}(\nbigm'),F_{\dagger}^j(\nbigm''),F_{\dagger}C)$
$(j\in\seisuu)$ on $Y$.
(See \cite[\S1.4]{Sabbah-pure-twistor}.)
When $\nbigt$ is equipped with an increasing filtration $W$,
let
$W_{j+\ell}(F_{\dagger}^j(\nbigt))$
denote the image of
$F_{\dagger}^j(W_{\ell}\nbigt)
\to F_{\dagger}^j(\nbigt)$
by which we obtain the filtration $W$
on $F_{\dagger}^j(\nbigt)$.

\subsubsection{Pure twistor $\nbigd$-modules}

A polarizable pure twistor $\nbigd$-module of weight $w$
on a complex manifold $X$
is an $\nbigr_X$-triple satisfying some conditions.
(See \cite{Sabbah-wild-twistor,Mochizuki-wild}.)
Let $\MT(X,w)$ denote the category of
polarizable pure twistor $\nbigd_X$-modules of weight $w$.

\begin{thm}[\cite{Sabbah-pure-twistor, mochi2, Mochizuki-wild}]
\label{thm;25.6.28.1}
Let $F:X\to Y$ be a projective morphism of complex manifolds.
For any $\nbigt\in\MT(X,w)$,
we have $F^j_{\dagger}(\nbigt)\in \MT(Y,w+j)$.
\hfill\qed
\end{thm}

There exists the full subcategory
$\MT_{\reg}(X,w)\subset \MT(X,w)$ of
regular polarizable pure twistor $\nbigd_X$-modules of weight $w$.
(See \cite{Sabbah-pure-twistor,mochi2} for
the regularity condition.)
We obtain the following generalization of 
Theorem \ref{thm;25.6.28.1}
in the regular case.
\begin{thm}[\cite{Mochizuki-L2}]
Let $F:X\to Y$ be a morphism of complex manifolds.
Let $\nbigt\in\MT_{\reg}(X,w)$.
Suppose that $X$ is K\"ahler,
and that the support of $\nbigt$ is proper over $Y$.
Then, we have $F^j_{\dagger}(\nbigt)\in \MT_{\reg}(Y,w+j)$.
\hfill\qed
\end{thm}

\subsubsection{Mixed twistor $\nbigd$-modules}

A mixed twistor $\nbigd$-module on $X$
is a filtered $\nbigr_X$-triple $(\nbigt,W)$
and satisfying some additional conditions.
(See \cite[\S7]{Mochizuki-MTM}.)
Let $\MTM(X)$ denote the category of
mixed twistor $\nbigd$-modules on $X$.
The following theorem is fundamental.
\begin{thm}[\cite{Mochizuki-MTM}]
\label{thm;25.6.28.2}
Let $F:X\to Y$ be a projective morphism of complex manifolds.
For $(\nbigt,W)\in\MTM(X)$,
we have $(F_{\dagger}^j(\nbigt),W)\in \MTM(Y)$ for any $j$.
\hfill\qed
\end{thm}

For any $(\nbigt,W)\in\MTM(X)$,
we have $\Gr^W_w(\nbigt)\in\MT(X,w)$.
There exists the full subcategory
$\MTM_{\reg}(X)\subset \MTM(X)$ of
mixed twistor $\nbigd_X$-modules $(\nbigt,W)$
such that
$\Gr^W_w(\nbigt)\in\MT_{\reg}(X,w)$.
Theorem \ref{thm;25.6.28.2}
is generalized as follows in the regular case.
\begin{thm}[\cite{Mochizuki-L2}]
Let $F:X\to Y$ be a morphism of complex manifolds.
Let $\nbigt\in\MTM_{\reg}(X)$.
Suppose that $X$ is K\"ahler,
and that the support of $\nbigt$ is proper over $Y$.
Then, we have $F^j_{\dagger}(\nbigt)\in \MTM_{\reg}(Y)$.
\hfill\qed
\end{thm}

\subsubsection{Pure and mixed twistor $\nbigd$-modules on a point}

Let $\pt$ denote the set of one point.
An $\nbigr_{\pt}$-module is an $\nbigo_{\cnum_{\lambda}}$-module.
Let $\sigma:\proj^1\to\proj^1$ be
the anti-holomorphic map defined by
$\sigma(\lambda)=(-\lambdabar)^{-1}$.
Let $\sigma:\cnum_{\lambda}^{\ast}\to\cnum_{\lambda}^{\ast}$
and $\sigma:\proj^1\setminus\{0\}\to\cnum_{\lambda}$
denote the induced maps.
Let $\nbigm',\nbigm''$ be locally free
$\nbigo_{\cnum_{\lambda}}$-modules of finite rank
with $\rank\nbigm'=\rank\nbigm''$.
Let 
\[
C:\nbigm'_{|\cnum_{\lambda}^{\ast}}
\otimes
 \sigma^{\ast}(\nbigm''_{|\cnum_{\lambda}^{\ast}})
 \to
 \nbigo_{\cnum_{\lambda}^{\ast}}
\]
be an $\nbigo_{\cnum_{\lambda}^{\ast}}$-homomorphism
which is perfect in the sense
the induced morphism
$\Psi_C:\sigma^{\ast}(\nbigm'')_{|\cnum_{\lambda}^{\ast}}
\to
(\nbigm')^{\lor}_{|\cnum_{\lambda}^{\ast}}$
is an isomorphism,
where
$(\nbigm')^{\lor}
=\nhom_{\nbigo_{\cnum_{\lambda}}}(\nbigm',\nbigo_{\cnum_{\lambda}})$.
Such a triple $(\nbigm',\nbigm'',C)$
is called a smooth $\nbigr_{\pt}$-triple.
For any smooth $\nbigr_{\pt}$-triple
$\nbigt=(\nbigm',\nbigm'',C)$,
we obtain
the locally free $\nbigo_{\proj^1}$-module
$\Upsilon(\nbigt)$
by gluing
$(\nbigm')^{\lor}$ on $\cnum_{\lambda}$
and 
$\sigma^{\ast}(\nbigm'')$ on $\proj^1\setminus\{0\}$
with $\Psi_C$.

A polarizable pure twistor $\nbigd$-module on $\pt$
is a smooth $\nbigr_{\pt}$-triple $\nbigt$
such that
$\Upsilon(\nbigt)$
is isomorphic to a direct sum of $\nbigo_{\proj^1}(w)$.
A mixed twistor $\nbigd$-module on $\pt$
is a smooth $\nbigr_{\pt}$-triple $\nbigt$
with a weight filtration $W$
such that
$\Gr^W_w(\nbigt)$
are pure of weight $w$.

Let $a_X:X\to\pt$ denote the canonical morphism.
We set
$\Omegatilde^k_{\nbigx/\cnum}:=
\lambda^{-k}(p_{\lambda}^{\ast}\Omega_X^k)$.
We have the exterior derivative
$d:\Omegatilde^k_{\nbigx/\cnum}\to\Omegatilde^{k+1}_{\nbigx/\cnum}$.
For a coherent $\nbigr_X$-module $\nbigm$ with compact support,
we obtain the complex of sheaves
$\nbigm\otimes \Omegatilde^{\bullet}_{\nbigx/\cnum}$
on $\nbigx$.
We have
\[
 a_{X\dagger}^{\ell}(\nbigm)=
 R^{\dim X+\ell}(\id\times a_X)_{\ast}\Bigl(
 \nbigm\otimes\Omegatilde^{\bullet}_{\nbigx/\cnum}
 \Bigr)
\]
as the $\nbigo_{\cnum_{\lambda}}$-module.

\begin{cor}
\label{cor;25.6.29.1}
Suppose that $\nbigm$ underlies
a regular mixed twistor $\nbigd$-module on $X$
with compact support. 
We also assume that $X$ is K\"ahler.
Then,
$a_{X\dagger}^{\ell}(\nbigm)$
are locally free $\nbigo_{\cnum_{\lambda}}$-modules.
\hfill\qed
\end{cor}

This is closely related with the $E_1$-degeneration property.
We consider the subcomplexes
$\nbigf_{\lambda}^j(\nbigm\otimes\Omegatilde^{\bullet}_{\nbigx/\cnum})
=\lambda^j\nbigm\otimes\Omegatilde^{\bullet}_{\nbigx/\cnum}$.
Because
$\lambda^j:a^{\ell}_{X\dagger}(\nbigm)
\to
 a^{\ell}_{X\dagger}(\nbigm)$
are monomorphisms,
the following is exact for any $\ell$ and $j$:
\begin{multline}
 0\lrarr
 R^{\dim X+\ell}(\id\times a_X)_{\ast}
 \nbigf_{\lambda}^j(\Omegatilde^{\bullet}_{\nbigx/\cnum}\otimes\nbigm)
 \lrarr
  R^{\dim X+\ell}(\id\times a_X)_{\ast}
 (\Omegatilde^{\bullet}_{\nbigx/\cnum}\otimes\nbigm)
 \lrarr
 \\
 R^{\dim X+\ell}(\id\times a_X)_{\ast}
 \bigl(
 (\Omegatilde^{\bullet}_{\nbigx/\cnum}\otimes\nbigm)/\nbigf_{\lambda}^j
 \bigr)
 \lrarr 0.
\end{multline}
Hence, the following is an epimorphism
for any $j$ and $\ell$:
\begin{equation}
 R^{\dim X+\ell}(\id\times a_X)_{\ast}
 \bigl(
 (\Omegatilde^{\bullet}_{\nbigx/\cnum}\otimes\nbigm)
 \big/\nbigf_{\lambda}^{j+1}
 \bigr)
 \lrarr
  R^{\dim X+\ell}(\id\times a_X)_{\ast}\bigl(
 (\Omegatilde^{\bullet}_{\nbigx/\cnum}\otimes \nbigm)/\nbigf_{\lambda}^j
 \bigr).
\end{equation}
This means the $E_1$-degeneration
of the spectral sequence associated with the filtration
$\nbigf_{\lambda}$
on $\Omegatilde^{\bullet}_{\nbigx/\cnum}\otimes\nbigm$.

\subsection{Integrable mixed twistor $\nbigd$-modules}

We set
$\nbigrtilde_X=\nbigr_X\langle\lambda^2\del_{\lambda}\rangle
\subset\nbigd_{\nbigx}$.
If $X$ is an open subset in $\cnum^n$,
we have
$\nbigrtilde_X=\nbigo_{\nbigx}
\langle \lambda\del_1,\ldots,\lambda\del_n,
\lambda^2\del_{\lambda}\rangle$.
By the identification $\vecS=\{e^{\sqrt{-1}\theta}\}$,
we obtain the vector field $\del_{\theta}$ on $\vecS$.

Let $\nbigm',\nbigm''$ be $\nbigrtilde_X$-modules.
Let $\nbigt=(\nbigm',\nbigm'',C)$
be an $\nbigr_X$-triple.
For any section $m'$ of $\nbigm'$
for $U\subset \vecS\times X$,
we set
$\del_{\theta}m'=\sqrt{-1}\lambda\del_{\lambda}m'
=(\sqrt{-1}\lambda\del_{\lambda}
-\sqrt{-1}\lambdabar\del_{\lambdabar})m'$.
Similarly, $\del_{\theta}m''$ is defined for
a section $m''$ of $\nbigm''$.
The $\nbigr$-triple $\nbigt$ is called integrable
if
\[
 \del_{\theta}C(m',\overline{\sigma^{-1}(m'')})
=C(\del_{\theta}m',\overline{\sigma^{-1}(m'')})
+C(m',\overline{\sigma^{-1}(\del_{\theta}m'')}).
\]
(See \cite[\S2.1.5]{Mochizuki-MTM} for integrable
$\nbigr_X$-triples,
which originally goes back to \cite{Sabbah-pure-twistor}.)
An integrable $\nbigr_X$-triple is called
$\nbigrtilde_X$-triple.
A morphism of $\nbigrtilde_X$-triples
$\nbigt_i=(\nbigm'_i,\nbigm''_i,C)$ $(i=1,2)$
is defined to be a morphism $(\varphi',\varphi'')$
of $\nbigr_X$-triples
such that $\varphi'$ and $\varphi''$
are $\nbigrtilde_X$-homomorphisms.
For a morphism of complex manifolds $F:X\to Y$
and for an $\nbigrtilde_X$-triple $\nbigt$,
the $\nbigr_Y$-triples
$F^j_{\dagger}(\nbigt)$
are naturally $\nbigrtilde_Y$-triples.

An integrable mixed twistor $\nbigd$-module on $X$
is a filtered $\nbigrtilde_X$-triple $(\nbigt,W)$
satisfying some conditions.
(See \cite[\S7.2.3]{Mochizuki-MTM}.)
Let $\MTM_{\reg}^{\integral}(X)$ denote
the category of
integrable mixed twistor $\nbigd_X$-modules
whose underlying mixed twistor $\nbigd_X$-modules are regular.
Let $\nbigc_{\reg}(X)$ denote
the full subcategory of $\nbigrtilde_X$-modules
underlying regular integrable mixed twistor $\nbigd_X$-modules,
i.e.,
an $\nbigrtilde_X$-module $\nbigm''$ is an object of $\nbigc_{\reg}(X)$
if and only if
there exists
$((\nbigm',\nbigm'',C),W)\in\MTM^{\integral}_{\reg}(X)$.

\subsection{$\nbigrtilde_X$-modules induced by Hodge modules}

Let $\nbigd_X$ denote the sheaf of
holomorphic linear differential operators on $X$.
Let $F_j(\nbigd_X)$ denote the subsheaf
of differential operators of degree at most $j$.
We set
$R^F(\nbigd_X):=\sum_{j\in\seisuu} \lambda^jF_j(\nbigd_X)$
and
$\Rtilde^F(\nbigd_X):=
R^F(\nbigd_X)\langle \lambda^2\del_{\lambda}\rangle$.

Let $M$ be a regular holonomic $\nbigd_X$-module.
Let $F(M)$ be a good filtration of $M$.
We obtain $R^F(\nbigd_X)$-module
$R^F(M)=\sum_{j\in\seisuu}\lambda^jF_j(M)$.
It is naturally an $\Rtilde^F(\nbigd_X)$-module.
By the analytification,
it induces an $\nbigrtilde_X$-module
denoted by $\nbigr^F(M)$.
In this way,
we obtain a functor
from the category of good filtered regular holonomic
$\nbigd_X$-modules
to the category of $\nbigrtilde_X$-modules.

\begin{lem}
If $(M,F)$ is a filtered regular holonomic $\nbigd$-module
underlying a mixed Hodge module,
we have $\nbigr^F(M)\in \nbigc_{\reg}(X)$.
(See {\rm\cite{saito1,saito2}} for Hodge modules.)
\end{lem}
\pf
There exists a natural functor
from the category of mixed Hodge modules on $X$
to $\MTM^{\integral}_{\reg}(X)$
as explained in \cite[\S13.5]{Mochizuki-MTM}.
In the level of filtered $\nbigd$-modules,
it is given as above.
\hfill\qed

\section{Main Theorem}

\subsection{Preliminary}

Let $X$ be a complex manifold.
Let $\nbigm\in\nbigc_{\reg}(X)$ be an $\nbigrtilde_X$-module
underlying an integrable regular mixed twistor $\nbigd$-module
induced by a mixed Hodge module.
Let $\Ch(\Xi_{\DR}\nbigm)\subset T^{\ast}X$
denote the characteristic variety
of the underlying $\nbigd_X$-module $\Xi_{\DR}(\nbigm)$.
Let $0_X:X\to T^{\ast}X$ denote the $0$-section.
We assume the following.
\begin{itemize}
 \item The set $\Cr(f)$ is compact.
 \item Any $c\neq 0$ is a regular value of $f$.
 \item $\Ch(\Xi_{\DR}(\nbigm))\cap df(X)\subset 0_X(\Cr(f))$.
\end{itemize}
The second condition implies $\Cr(f)\subset f^{-1}(0)$.
Because the characteristic varieties are cone,
the third condition implies
$\Ch(\Xi_{\DR}(\nbigm))\cap (\alpha df)(X)\subset
0_X(\Cr(f))$
for any non-zero constant $\alpha$.

\subsubsection{Cohomology group of the restriction to $\lambda=0$}

Let $\nbigl(f)$ denote the $\nbigrtilde_X$-module
given by
$\nbigo_{\nbigx}$
with the meromorphic integrable connection
$d+d(\lambda^{-1}f)$.
We obtain the $\nbigrtilde_X$-module
$\nbigm_f=\nbigm\otimes\nbigl(f)$ on $\nbigx$,
and the complex of sheaves
$\nbigm_f\otimes\Omegatilde^{\bullet}_{\nbigx/\cnum}$
on $\nbigx$.
We recall that
$\nbigm_f\in\nbigc_{\reg}(X)$.
Let $\iota_{\lambda}:\{0\}\times X\to \cnum_{\lambda}\times X$
denote the inclusion.
We obtain the complex
of coherent $\Sym \Theta_X$-modules
$\iota_{\lambda}^{\ast}(\nbigm_f\otimes\Omegatilde^{\bullet}_{\nbigx/\cnum})$.
It induces a complex of
coherent $\nbigo_{T^{\ast}X}$-modules
denoted by
$\bigl(
\iota_{\lambda}^{\ast}(\nbigm_f\otimes\Omegatilde^{\bullet}_{\nbigx/\cnum})
\bigr)^{\sim}$.

\begin{lem}
\label{lem;24.6.25.10}
The cohomological support of 
$\bigl(
\iota_{\lambda}^{\ast}(\nbigm_f\otimes\Omegatilde^{\bullet}_{\nbigx/\cnum})
\bigr)^{\sim}$
is contained in $0_X(\Cr(f))$.
As a result,
\[
H^{\ast}\bigl(
X,\iota_{\lambda}^{\ast}(\nbigm_f\otimes\Omegatilde^{\bullet}_{\nbigx/\cnum})
\bigr) 
\]
are finite dimensional.
\end{lem}
\pf
Let $(\iota_{\lambda}^{\ast}\nbigm)^{\sim}$
and $(\iota_{\lambda}^{\ast}\nbigl(f))^{\sim}$
denote the coherent $\nbigo_{T^{\ast}X}$-modules
induced by
$\iota_{\lambda}^{\ast}\nbigm$
and $\iota_{\lambda}^{\ast}\nbigl(f)$,
respectively.
The support of
$(\iota_{\lambda}^{\ast}\nbigm)^{\sim}$
is the characteristic variety $\Ch(\Xi_{\DR}\nbigm)$.
The support of 
$(\iota_{\lambda}^{\ast}\nbigl(f))^{\sim}$
is the image of $df(X)$.
The support of
$(\iota_{\lambda}^{\ast}\nbigm_f)^{\sim}$
is $df(X)+\Ch(\Xi_{\DR}(\nbigm))$ in $T^{\ast}X$.

Let $\omega_X$ denote the canonical bundle of $X$.
Let 
$\Sym\Theta_X
 \otimes \bigwedge^{\bullet}\Theta_X
 \otimes\omega_X$
be the Koszul resolution of $\nbigo_X$
by $\Sym\Theta_X$-free modules.
It induces an $\nbigo_{T^{\ast}X}$-free resolution
$\nbigo_{T^{\ast}X}\otimes
\pi^{\ast}\bigl(
\bigwedge^{\bullet}\Theta_X\otimes\omega_X\bigr)$
of $0_{X\ast}(\nbigo_X)$,
where $\pi:T^{\ast}X\to X$ denotes the projection.

We have
\[
\bigl(
\iota_{\lambda}^{\ast}(\nbigm_f\otimes\Omegatilde^{\bullet}_{\nbigx/\cnum})
\bigr)^{\sim}
\simeq
(\iota_{\lambda}^{\ast}\nbigm_f)^{\sim}
\otimes_{\nbigo_{T^{\ast}X}}
\Bigl(
\nbigo_{T^{\ast}X}\otimes
\pi^{\ast}\Bigl(
\bigwedge^{\bullet} \Theta_X\otimes\omega_X
\Bigr)
\Bigr).
\]
Hence, the cohomological support is contained in
the intersection of
the support of 
$(\iota_{\lambda}^{\ast}\nbigm_f)^{\sim}$
and the $0_X(X)$,
which is contained in $0_X(\Cr(f))$.
\hfill\qed

\subsubsection{Cohomology group of the vanishing cycle sheaf}

Let $\iota_f:X\to X\times\cnum_t$
be the graph embedding,
i.e.,
$\iota_f(x)=(x,f(x))$.
There exists the $V$-filtration
$V(\iota_{f\dagger}(\Xi_{\DR}(\nbigm)))$
along $t$.
We obtain the regular holonomic $\nbigd_X$-module
\[
 \phi_f(\Xi_{\DR}(\nbigm))
:=\bigoplus_{-1< a\leq 0}
 \Gr^V_a\bigl(
 \iota_{f\dagger}\Xi_{\DR}(\nbigm)
 \bigr).
\]

\begin{lem}
\label{lem;25.8.26.10}
The support of $\phi_f(\Xi_{\DR}(\nbigm))$
is contained in the compact subset $\Cr(f)$.
As a result,
the cohomology group
 $H^{\ast}\bigl(
 X,\phi_f(\Xi_{\DR}\nbigm)\otimes\Omega^{\bullet}_X
 \bigr)$
is finite dimensional.
\end{lem}
\pf
The third condition implies that
$\Xi_{\DR}(\nbigm)$ is non-characteristic to
the hypersurfaces $f^{-1}(c)$ on $X\setminus\Cr(f)$.
Hence, $\phi_f(\Xi_{\DR}(\nbigm))=0$ on $X\setminus \Cr(f)$.
\hfill\qed

\subsubsection{$E_1$-degeneration and long exact sequences}
\label{subsection;25.8.25.10}

We recall that the $E_1$-degeneration condition
of a spectral sequence associated with a filtered complex
for the convenience of readers.
Let $\nbiga$ be an abelian category.
Let $(K^{\bullet},d)$ be a complex
with a decreasing filtration $F^{\bullet}(K^{\bullet})$.
To simplify the notation,
we set
$F^{p,p+r}(K^{j}):=F^p(K^{j})/F^{p+r}(K^{j})$.

\begin{lem}
Let $r_0\in\seisuu_{\geq 1}$.
The following conditions are equivalent.
\begin{description}
 \item[$A(r_0)$:] $H^j(F^{p,p+r}(K^{\bullet}))
       \to
       H^j(F^{p,p+1}(K^{\bullet}))$
    are epimorphisms for any $j,p\in\seisuu$
    and $0\leq r\leq r_0$.
 \item[$B(r_0)$:] $H^j(F^{p,p+r}(K^{\bullet}))
       \to
       H^j(F^{p,p+r-1}(K^{\bullet}))$
       are epimorphisms for any $j,p\in\seisuu$
      and $0\leq r\leq r_0$.
\end{description}
\end{lem}
\pf
It is easy to see that $B(r_0)$ implies $A(r_0)$.
Suppose that $A(r_0)$ holds.
We shall prove
$H^j(F^{p,p+r}(K^{\bullet}))
\to
H^j(F^{p,p+r-1}(K^{\bullet}))$
are epimorphisms for any $j,p\in\seisuu$
and $0\leq r\leq r_0$,
by an induction on $r$.
We have the following commutative diagram:
\[
\begin{CD}
 @. H^j(F^{p,p+r-1})
 @>>> H^j(F^{p,p+1}) \\
 @. @V{a}VV @V{b}VV \\
 0 @>>>
 H^{j+1}(F^{p+r-1,p+r})
 @>{c}>>
 H^{j+1}(F^{p+1,p+r})
 @>>>
 H^{j+1}(F^{p+1,p+r-1})
 @>>> 0.
\end{CD}
\]
By $A(r_0)$,
we have $b=0$.
By the assumption of the induction on $r$,
$c$ is a monomorphism.
Hence, we obtain $a=0$,
and the induction can proceed.
\hfill\qed

\vspace{.1in}

Recall that the spectral sequence
$E_r^{p,q}$
for the filtered complex
$F^{\bullet}(K^{\bullet})$
is given as
\[
 E_r^{p,q}
 =Z_r^{p,q}\big/
 (Z_{r-1}^{p+1}+dZ^{p-r+1,q+r-2}_{r-1})
\]
by setting
$Z_r^{p,q}=
 \Ker\Bigl(
 d:F^{p}K^{p+q}\lrarr
 F^{p,p+r}(K^{p+q+1})
 \Bigr)$.
(See \cite{Gelfand-Manin} for more details.)
There exist the natural morphisms
$d_r:E_r^{p,q}\to E_{r}^{p+r,q-r+1}$
such that $d_r\circ d_r=0$,
induced by $d:F^pK^{p+q}\to F^pK^{p+q+1}$.
There exist natural isomorphisms
\[
 E_{r+1}^{p,q}
 \simeq
 \Ker\Bigl(
 E_r^{p,q}\to E_r^{p+r,q-r+1}
 \Bigr)\Big/
 \Image\Bigl(
 E_r^{p-r,q+r-1}
 \to
 E_r^{p,q}
 \Bigr).
\]

\begin{lem}
We have $d_r=0$ for any $1\leq r\leq r_0-1$
if and only if $A(r_0)$ holds.
\end{lem}
\pf
We shall use an induction on $r_0$.
Suppose that
$d_r=0$ $(1\leq r\leq r_0-2)$
and
that $A(r_0-1)$ holds.
We have
$E_r^{p,q}=E_1^{p,q}
=H^{p+q}(F^{p,p+1}(K^{\bullet}))$ for $1\leq r\leq r_0-1$.
We consider
\[
 \begin{CD}
 @. @. H^j(F^{p,p+1})\\
 @. @. @V{a}VV \\
0 @>>>
 H^{j+1}(F^{p+r_0-1,p+r_0})
 @>{c}>>  
 H^{j+1}(F^{p+1,p+r_0})
 @>{b}>>
  H^{j+1}(F^{p+1,p+r_0-1})
  @>>> 0.
\end{CD}
\]
By $A(r_0-1)$, $c$ is a monomorphism,
and we have $b\circ a=0$.
There exists a unique morphism
$\varphi:H^j(F^{p,p+1})\to H^{j+1}(F^{p+r_0-1,p+r_0})$
such that
$c\circ\varphi=b$.
By the construction,
$\varphi=d_{r_0-1}$.
Hence, $A(r_0)$ holds if and only if $d_{r_0-1}=0$.
\hfill\qed

\vspace{.1in}
We say that the spectral sequence degenerates at the $E_1$-level
if $d_r=0$ for any $r\geq 1$.
We obtain the following proposition.
\begin{prop}
The $E_1$-degeneration holds if and only if
one of the following equivalent conditions holds. 
\begin{itemize}
 \item $H^j(F^{p,p+r}(K^{\bullet}))
       \to
       H^j(F^{p,p+r-1}(K^{\bullet}))$
       are epimorphisms for any $j,p\in\seisuu$
      and $1\leq r$.
\hfill\qed
\end{itemize}
\end{prop}

\subsection{Refinement of the Theorem \ref{thm;23.3.29.1}}

\subsubsection{Main Theorem}

Let us explain a refined statement
of Theorem \ref{thm;23.3.29.1}.
We define the filtration
$\nbigf_{\lambda}^k(\nbigm\otimes\Omegatilde^{\bullet}_{\nbigx/\cnum})
=\lambda^k
(\nbigm\otimes\Omegatilde^{\bullet}_{\nbigx/\cnum})$
for $k\in\seisuu_{\geq 0}$.
Note that
$\Gr_{\nbigf_{\lambda}}^k\bigl(
 \nbigm\otimes\Omegatilde^{\bullet}_{\nbigx/\cnum}
\bigr)$
is isomorphic to 
$\iota_{\lambda}^{\ast}(\nbigm_f\otimes\Omegatilde^{\bullet}_{\nbigx/\cnum})$
for any $k\geq 0$.
Let us state the main theorem of this paper.

\begin{thm}
\label{thm;23.3.29.2}
The $E_1$-degeneration holds
for the filtered complex
$\nbigf_{\lambda}(\nbigm\otimes\Omegatilde^{\bullet}_{\nbigx/\cnum})$
with respect to the push-forward by
the projection $\nbigx\to\cnum_{\lambda}$.
In other words,
we obtain the following exact sequences
for any $\ell$  and $j$:
\begin{multline}
\label{eq;25.6.30.1}
 0\to
 H^j\Bigl(
 X,
 \Gr^{\nbigf_{\lambda}}_{\ell}
 \bigl(
\nbigm_f\otimes\Omegatilde^{\bullet}_{\nbigx/\cnum}
\bigr)
 \Bigr)
 \to
  H^j\Bigl(
 X,
 \bigl(
\nbigm_f\otimes\Omegatilde^{\bullet}_{\nbigx/\cnum}
 \bigr)\big/
 \nbigf_{\lambda}^{\ell}\bigl(
\nbigm_f\otimes\Omegatilde^{\bullet}_{\nbigx/\cnum}
 \bigr)
 \Bigr)
\\
 \to
  H^j\Bigl(
 X,
 \bigl(
\nbigm_f\otimes\Omegatilde^{\bullet}_{\nbigx/\cnum}
 \bigr)\big/
 \nbigf_{\lambda}^{\ell+1}\bigl(
\nbigm_f\otimes\Omegatilde^{\bullet}_{\nbigx/\cnum}
 \bigr)
 \Bigr)
\to 0.
\end{multline} 
\end{thm}
We obtain Theorem \ref{thm;23.3.29.1}
from Theorem \ref{thm;23.3.29.2}
as the special case $\nbigm=\nbigo_{\nbigx}$.
We shall also prove the following.
\begin{thm}
We have
$\dim H^j\bigl(X,
\iota_{\lambda}^{\ast}(\nbigm_f\otimes\Omegatilde^{\bullet}_{\nbigx/\cnum})
\bigr)
=\dim H^{j}\bigl(
X,\phi_f(\Xi_{\DR}\nbigm)\otimes
\Omega^{\bullet}_X
\bigr)$
for any $j$.
\end{thm}
\begin{cor}
In particular,
we have 
$\dim H^j\bigl(X,
 (\Omega^{\bullet}_X,df)
 \bigr)
 =\dim H^{j}(X,\Omega^{\bullet}_X\otimes
 \phi_f(\nbigo_X))$.
\end{cor}

\begin{cor}
There exists an isomorphism
\begin{equation}
\varprojlim_{\ell}
H^j\Bigl(
 X,
\bigl(
\nbigm_f\otimes\Omegatilde^{\bullet}_{\nbigx/\cnum}
\bigr)
\Big/
\nbigf^{\ell}_{\lambda}
\bigl(
\nbigm_f\otimes\Omegatilde^{\bullet}_{\nbigx/\cnum}
\bigr)
 \Bigr)
 \simeq
  H^j\Bigl(
 X,
 \iota_{\lambda}^{\ast}\bigl(
 \nbigm_f\otimes
 \Omegatilde^{\bullet}_{\nbigx/\cnum}
 \bigr)
 \Bigr)
\otimes
 \cnum[\![\lambda]\!].
\end{equation}
\hfill\qed
\end{cor}

\subsubsection{Completion}

We naturally regard
$\bigl(
\nbigm_f
\bigl/
\nbigf^j_{\lambda}\nbigm_f
\bigr)
\otimes\Omegatilde^{\bullet}_{\nbigx/\cnum}$
as the complexes of sheaves on $X$.
We obtain the following complex of sheaves on $X$:
\[
\widehat{\nbigm_f\otimes\Omegatilde^{\bullet}_{\nbigx/\cnum}}=
\varprojlim_j
\bigl(
\nbigm_f
\big/
\nbigf^j_{\lambda}\nbigm_f
\bigr)
\otimes\Omegatilde^{\bullet}_{\nbigx/\cnum}.
\]
\begin{cor}
\label{cor;25.7.2.1}
There exists an isomorphism
\[
H^j\Bigl(
 X,
\widehat{\nbigm_f\otimes\Omegatilde^{\bullet}_{\nbigx/\cnum}}
\Bigr)
\simeq
 H^j\Bigl(
 X,
 \iota_{\lambda}^{\ast}\bigl(
 \nbigm_f\otimes
 \Omegatilde^{\bullet}_{\nbigx/\cnum}
 \bigr)
 \Bigr)
\otimes
 \cnum[\![\lambda]\!].
\]
In particular,
$H^j\Bigl(
 X,
\widehat{\nbigm_f\otimes\Omegatilde^{\bullet}_{\nbigx/\cnum}}
\Bigr)$
are free $\cnum[\![\lambda]\!]$-modules of finite rank.
\end{cor}
\pf
Let $\Omega^{0,q}_X$
denote the sheaf of smooth $(0,q)$-forms on $X$.
We obtain the following double complex
\[
\bigl(
 \widehat{
 \nbigm_f\otimes\Omegatilde^{\bullet}_{\nbigx/\cnum}}
\bigr)
 \otimeshat_{\nbigo_X}
 \Omega^{0,\bullet}_X
 :=
\varprojlim_j
\left(
\bigl(
\nbigm_f\big/
\nbigf^j_{\lambda}\nbigm_f
\bigr)
\otimes\Omegatilde^{\bullet}_{\nbigx/\cnum}
\otimes_{\nbigo_X}
\Omega_X^{0,\bullet}
\right).
\]
Let
$\Tot\Bigl(
\bigl(
\widehat{
 \nbigm_f\otimes\Omegatilde^{\bullet}_{\nbigx/\cnum}}
\bigr)
 \otimeshat_{\nbigo_X}
 \Omega^{0,\bullet}_X\Bigr)$
 denote the total complex.
 \begin{lem}
\label{lem;25.8.25.2}
The natural morphism is a fine resolution:
\begin{equation}
\label{eq;25.8.25.1}
 \widehat{
 \nbigm_f\otimes\Omegatilde^{\bullet}_{\nbigx/\cnum}}
 \lrarr
 \Tot\Bigl(
\bigl(
 \widehat{
 \nbigm_f\otimes\Omegatilde^{\bullet}_{\nbigx/\cnum}}
\bigr)
 \otimeshat_{\nbigo_X}
 \Omega^{0,\bullet}_X\Bigr).
\end{equation}
 \end{lem}
\pf
It is enough to prove that
(\ref{eq;25.8.25.1}) is a quasi-isomorphism
locally around any point of $X$.
Recall that the sheaf of $C^{\infty}$-functions on $X$
is flat over $\nbigo_X$
according to \cite{malgrange2}.

Let $G$ be any pseudo-coherent $\nbigo_X$-module.
(See \cite[Appendix A]{kashiwara_text}
for pseudo-coherent sheaves.)
We have
$G\otimes^L_{\nbigo_X}\Omega^{0,q}_X
\simeq
G\otimes_{\nbigo_X}\Omega^{0,q}_X$.
Let $F_{\bullet}\to G$ be a free resolution of $G$.
The natural morphisms
$F_{\bullet}\to
\Tot(F_{\bullet}\otimes_{\nbigo_X}\Omega^{0,\bullet}_X)
\to G\otimes_{\nbigo_X}\Omega^{0,\bullet}_X$
are quasi-isomorphisms.
Hence, 
$G\to G\otimes_{\nbigo_X}\Omega^{0,\bullet}_X$
is a quasi-isomorphism.

Let $\nbigg$ be a pseudo-coherent $\nbigo_{\nbigx}$-module
flat over $\nbigo_{\cnum_{\lambda}}$.
We naturally regard $\nbigg/\lambda^j\nbigg$
as $\nbigo_X$-modules.
Let $\pi_{j+1}:\nbigg/\lambda^{j+1}\nbigg\to\nbigg/\lambda^j\nbigg$
denote the projections.
We obtain the quasi-isomorphisms
\[
 \nbigg/\lambda^j\nbigg
 \to
 (\nbigg/\lambda^j\nbigg)\otimes_{\nbigo_X}\Omega_X^{0,\bullet}.
\]
For any open subset $U\subset X$,
the morphisms
$H^0\bigl(U,(\nbigg/\lambda^{j+1}\nbigg)\otimes\Omega^{0,q}_X\bigr)
 \to
 H^0\bigl(U,(\nbigg/\lambda^{j}\nbigg)\otimes\Omega^{0,q}_X\bigr)$
are surjective
because $\iota_{\lambda}^{\ast}(\nbigg)\otimes\Omega^{0,q}_X$ is fine.
Let $a_j^{q+1}\in
H^0\bigl(U,(\nbigg/\lambda^j\nbigg)\otimes\Omega^{0,q+1}\bigr)$ $
(j=1,2,\ldots)$
be sections such that
$\delbar(a_j^{q+1})=0$
and $\pi_{j+1}(a_{j+1}^{q+1})=a_j^{q+1}$.
Let us construct
$b_j^q\in H^0\bigl(U,(\nbigg/\lambda^j)\otimes\Omega^{0,q}_X\bigr)$
$(j=1,2,\ldots)$
such that
$\delbar b_j^q=a_j^q$
and $\pi_{j+1}(b_{j+1}^q)=b_j^q$
inductively on $j$.
Suppose that we have already constructed $b_{j}^q$.
There exists
$c_{j+1}^q\in H^0\bigl(U,(\nbigg/\lambda^{j+1}\nbigg)\otimes\Omega^{0,q}_X
\bigr)$
such that
$\pi_{j+1}(c_{j+1}^q)=b_j^q$.
We obtain
\[
d_{j+1}^{q+1}=a_{j+1}^{q+1}-\delbar c_{j+1}^q
\in
H^0\bigl(U,
 (\lambda^j\nbigg/\lambda^{j+1}\nbigg)
 \otimes\Omega^{0,q+1}
\bigr)
\simeq H^0\bigl(U,
 \iota_{\lambda}^{\ast}\nbigg
 \otimes\Omega^{0,q+1}
\bigr)
\]
such that
$\delbar(d_{j+1}^{q+1})=0$.
There exists
$e_{j+1}^q\in
H^0\bigl(U,(\lambda^j\nbigg/\lambda^{j+1}\nbigg)\otimes
\Omega^{0,q}_X\bigr)$
such that
$\delbar e_{j+1}^q=d_{j+1}^{q+1}$.
By setting $b_{j+1}^q=c_{j+1}^q+e_{j+1}^q$,
the induction can proceed.

As a result,
the natural morphism 
\[
 \varprojlim_{j}\bigl(
 \nbigg/\lambda^j\nbigg
 \bigr)
 \to
 \varprojlim_{j}
 \Bigl(
 \bigl(
 \nbigg/\lambda^j\nbigg
 \bigr)
 \otimes\Omega^{0,\bullet}_X
 \Bigr)
\]
is a quasi-isomorphism.
Then, we obtain the claim of Lemma \ref{lem;25.8.25.2}.
\hfill\qed

\vspace{.1in}
We obtain Corollary \ref{cor;25.7.2.1}
from the following lemma.
\begin{lem}
The natural morphisms
\[
 H^j\Bigl(X, \widehat{
 \nbigm_f\otimes\Omegatilde^{\bullet}_{\nbigx/\cnum}}
  \Bigr)
 \to
  \varprojlim_{\ell}
 H^j\Bigl(
 X,
\bigl(
 \nbigm_f\big/
 \nbigf^{\ell}_{\lambda}\nbigm_f
 \bigr)
 \otimes\Omegatilde^{\bullet}_{\nbigx/\cnum}
 \Bigr)
\]
are isomorphisms.
\end{lem}
\pf
Let $\pi_{\ell+1}$ denote the projection
induced by
$\nbigm_f/\nbigf_{\lambda}^{\ell+1}\nbigm_f
\to
\nbigm_f/\nbigf_{\lambda}^{\ell}\nbigm_f$.
Let $d$ denote the differential of the complex
$\Tot^q\Bigl(
\bigl(
\nbigm_f
\big/
\nbigf^{\ell}_{\lambda}\nbigm_f
\bigr)
\otimes\Omegatilde^{\bullet}_{\nbigx/\cnum}
\otimes\Omega^{0,\bullet}_X
\Bigr)$.

Let us consider cohomology classes
\[
\alpha_{\ell}\in
 H^q\Bigl(
 X,
\bigl(
\nbigm_f\big/
\nbigf^{\ell}_{\lambda}\nbigm_f
\bigr)
\otimes\Omegatilde^{\bullet}_{\nbigx/\cnum}
\Bigr)
\quad
(\ell=1,2,\ldots)
\]
such that 
$\pi_{\ell+1}(\alpha_{\ell+1})=\alpha_{\ell}$.
We shall construct cocycles
\[
a_{\ell}
\in 
  H^0\left(
  X,
\Tot^q\Bigl(
\bigl(
\nbigm_f
\big/
\nbigf^{\ell}_{\lambda}\nbigm_f
\bigr)
\otimes\Omegatilde^{\bullet}_{\nbigx/\cnum}
\otimes\Omega^{0,\bullet}_X
\Bigr)
\right)
\quad
(\ell=1,2,\ldots)
\]
such that
$a_{\ell}$ are representatives of $\alpha_{\ell}$
and that $\pi_{\ell+1}(a_{\ell+1})=a_{\ell}$
inductively on $\ell$.
Suppose that we have already constructed $a_{\ell}$.
There exists 
$a'_{\ell+1}
 \in
  H^0\left(
  X,
\Tot^q\Bigl(
\bigl(
\nbigm_f\big/
\nbigf^{\ell+1}_{\lambda}\nbigm_f
\bigr)
\otimes\Omegatilde^{\bullet}_{\nbigx/\cnum}
\otimes\Omega^{0,\bullet}_X
\Bigr)
\right)$
such that
$\pi_{\ell+1}(a'_{\ell+1})=a_{\ell}$.
We obtain a cocycle
\[
 d(a'_{\ell+1})
\in
H^0\left(
X,
\Tot^{q+1}\Bigl(
\bigl(
\nbigf^{\ell}_{\lambda}\nbigm_f
\big/
\nbigf^{\ell+1}_{\lambda}\nbigm_f
\bigr)
\otimes\Omegatilde^{\bullet}_{\nbigx/\cnum}
\otimes\Omega^{0,\bullet}_X
\Bigr)
\right).
\]
Because $\pi_{\ell+1}(\alpha_{\ell+1})=\alpha_{\ell}$,
there exists
\[
 b_{\ell+1}\in
  H^0\left(
  X,
\Tot^q\Bigl(
\bigl(
\nbigf^{\ell}_{\lambda}\nbigm_f
\big/
\nbigf^{\ell+1}_{\lambda}\nbigm_f
\bigr)
\otimes\Omegatilde^{\bullet}_{\nbigx/\cnum}
\otimes\Omega^{0,\bullet}_X
\Bigr)
\right)
\]
such that
$d(b_{\ell+1})=d(a'_{\ell+1})$.
Let $\alpha'_{\ell+1}$
denote the cohomology class of $a'_{\ell+1}-b_{\ell+1}$.
Because
$\alpha_{\ell+1}-\alpha'_{\ell+1}$
comes from
$H^q\Bigl(
 X,\Tot\Bigl(
 (\nbigf_{\lambda}^{\ell}\nbigm_f
 \big/
 \nbigf_{\lambda}^{\ell+1}\nbigm_f)
 \otimes\Omegatilde_{\nbigx/\cnum}^{\bullet}\otimes\Omega_X^{0,\bullet}
 \Bigr)
 \Bigr)$,
there exists $a_{\ell+1}$
with the desired property.

Let us consider coboundaries
\[
 a_{\ell}
 \in
 H^0\Bigl(
 X,\Tot^{q}\Bigl(
 \bigl(
 \nbigm_f/
 \nbigf_{\lambda}^{\ell}\nbigm_f
 \bigr)
 \otimes\Omegatilde^{\bullet}_X
 \otimes\Omega_X^{0,\bullet}
 \Bigr)
 \Bigr)
 \quad
 (\ell=1,2,\ldots)
\]
such that
$\pi_{\ell+1}(a_{\ell+1})=a_{\ell}$.
Let us construct
\[
 b_{\ell}
 \in
 H^0\Bigl(
 X,\Tot^{q-1}\Bigl(
\bigl(
 \nbigm_f
\big/
 \nbigf_{\lambda}^{\ell}\nbigm_f
 \bigr)
 \otimes\Omegatilde^{\bullet}_X
 \otimes\Omega_X^{0,\bullet}
 \Bigr)
 \Bigr)
  \quad
 (\ell=1,2,\ldots)
\]
such that
$d(b_{\ell})=a_{\ell}$
and
$\pi_{\ell+1}(b_{\ell+1})= b_{\ell}$
inductively on $\ell$.
Suppose that we have already constructed $b_{\ell}$.
There exists
$b'_{\ell+1}
 \in
 H^0\Bigl(
 X,\Tot^{q-1}\Bigl(
 \bigl(
 \nbigm_f
 \big/
 \nbigf_{\lambda}^{\ell+1}\nbigm_f
 \bigr)
 \otimes\Omegatilde^{\bullet}_X
 \otimes\Omega_X^{0,\bullet}
 \Bigr)
 \Bigr)$
such that $\pi_{\ell+1}(b'_{\ell+1})= b_{\ell}$.
We obtain a cocycle
\[
 a_{\ell+1}-d(b'_{\ell+1})
 \in
 H^0\Bigl(
 X,\Tot^{q}\Bigl(
\bigl(
 \nbigf_{\lambda}^{\ell}\nbigm_f\big/
 \nbigf_{\lambda}^{\ell+1}\nbigm_f
\bigr)
 \otimes\Omegatilde^{\bullet}_X
 \otimes\Omega_X^{0,\bullet}
 \Bigr)
 \Bigr).
\]
By Theorem \ref{thm;23.3.29.2}, the natural morphism
\[
 H^q\bigl(
 X,\Gr^{\ell}_{\nbigf_{\lambda}}\nbigm
 \otimes\Omegatilde_{\nbigx/\cnum}^{\bullet}
 \bigr)
 \lrarr
 H^q\bigl(
 X,(\nbigm/\nbigf_{\lambda}^{\ell+1}\nbigm)
 \otimes\Omegatilde_{\nbigx/\cnum}^{\bullet}
 \bigr)
\]
is injective.
Hence, there exists 
\[
 c_{\ell+1}
 \in
 H^0\Bigl(
 X,\Tot^{q-1}\Bigl(
\bigl(
 \nbigf_{\lambda}^{\ell}\nbigm_f\big/
 \nbigf_{\lambda}^{\ell+1}\nbigm_f
\bigr)
 \otimes\Omegatilde^{\bullet}_X
 \otimes\Omega_X^{0,\bullet}
 \Bigr)
 \Bigr)
\]
such that
$d(c_{\ell+1})=a_{\ell+1}-d(b'_{\ell+1})$.
By setting
$b_{\ell+1}=b'_{\ell+1}+c_{\ell+1}$,
the induction can proceed.
\hfill\qed

\subsubsection{Reformulation}

Let $\iota_f:X\to X\times\cnum_t$ denote the graph embedding.
We obtain
the $\nbigrtilde_{X\times\cnum_t}$-module
$\nbigmtilde=\iota_{f\dagger}\nbigm$ on $X\times\cnum_t$.
The following lemma is standard.

\begin{lem}
There exists a natural isomorphism
$\nbigmtilde\otimes\nbigl(t)
 \simeq
 \iota_{f\dagger}(\nbigm_f)$. 
As a result,
there exists a natural quasi-isomorphism:
\[
 \nbigmtilde\otimes\nbigl(t)
 \otimes\Omegatilde^{\bullet}
 _{\cnum_{\lambda}\times X\times\cnum_t/\cnum_{\lambda}}
 \simeq
 \iota_{f\ast}\bigl(
 \nbigm_f\otimes\Omegatilde^{\bullet}_{\nbigx/\cnum_{\lambda}}
 \bigr)[-1].
\]
It induces the quasi-isomorphisms of the subcomplexes
\[
 \nbigf^k_{\lambda}\bigl(
 \nbigmtilde\otimes\nbigl(t)
 \otimes\Omegatilde^{\bullet}
 _{\cnum_{\lambda}\times X\times\cnum_t/\cnum_{\lambda}}
 \bigr)
 \simeq
 \iota_{f\ast}\bigl(
 \nbigf^k_{\lambda}\bigl(
 \nbigm_f\otimes\Omegatilde^{\bullet}_{\nbigx/\cnum_{\lambda}}
 \bigr)
 \bigr)[-1]
\]
and the quotient complexes:
\[
 \Gr^{k}_{\nbigf_{\lambda}}\bigl(
 \nbigmtilde\otimes\nbigl(t)
 \otimes\Omegatilde^{\bullet}
 _{\cnum_{\lambda}\times X\times\cnum_t/\cnum_{\lambda}}
 \bigr)
 \simeq
 \iota_{f\ast}\bigl(
 \Gr^k_{\nbigf_{\lambda}}\bigl(
 \nbigm_f\otimes\Omegatilde^{\bullet}_{\nbigx/\cnum_{\lambda}}
 \bigr)
 \bigr)[-1].
\] 
\hfill\qed
\end{lem}

Let $\pi_{01}:
\cnum_{\lambda}\times X\times\cnum_t\to \cnum_{\lambda}\times X$
denote the projection.
Let $V\nbigr_{X\times\cnum_t}\subset \nbigr_{X\times\cnum_t}$
denote
the sheaf of subalgebras
generated by
$\pi_{01}^{\ast}\nbigr_X$ and $t\deldel_t$.
Because
$\nbigmtilde\in\nbigc_{\reg}(X\times\cnum_t)$,
$\nbigmtilde$ has a $V$-filtration along $t$,
that is an increasing and exhaustive filtration
$V_{a}(\nbigmtilde)$ $(a\in\real)$ of $\nbigmtilde$
by coherent $V\nbigr_{X\times\cnum_t}$-submodules
satisfying the following conditions.
\begin{itemize}
 \item For any $a\in\real$,
       there exists $\epsilon>0$
       such that $V_a(\nbigmtilde)=V_{a+\epsilon}(\nbigmtilde)$.
 \item $\Gr^V_a=V_a/V_{<a}$
       are flat over $\nbigo_{\cnum_{\lambda}}$.
 \item $tV_a\subset V_{a-1}$,
       and $tV_a=V_{a-1}$ if $a<0$.
 \item $\deldel_t V_a\subset V_{a+1}$,
       and the induced morphisms
       $\deldel_t:\Gr^V_a\lrarr \Gr^V_{a+1}$
       are isomorphisms if $a>-1$.
 \item The induced actions of $-\deldel_tt-\lambda a$
       on $\Gr^V_a(\nbigmtilde)$
       are locally nilpotent.
\end{itemize}
Because $\nbigmtilde$ is induced by a mixed Hodge module,
$V_a(\nbigmtilde)$ are coherent over $\pi_{01}^{\ast}\nbigr_X$.

We obtain the following complex:
\[
 V_{-1}(\nbigmtilde)
 \lrarr
 \lambda^{-1}V_0(\nbigmtilde)\otimes dt,
 \quad
 s\longmapsto (\del_t+\lambda^{-1})s\,dt.
\]
The first term sits in the degree $0$.
It extends to the following double complex:
\begin{equation}
\label{eq;25.6.25.1}
  V_{-1}(\nbigmtilde)
  \otimes
  \pi_{01}^{\ast}\Omegatilde^{\bullet}_{\nbigx/\cnum}
 \lrarr
 \lambda^{-1}V_0(\nbigmtilde)\otimes dt
 \otimes
 \pi_{01}^{\ast}\Omegatilde^{\bullet}_{\nbigx/\cnum}.
\end{equation}
Let $\gbigc_1(\nbigm)$ denote the complex of sheaves
on $\cnum_{\lambda}\times X\times\cnum_t$
obtained as the total complex of (\ref{eq;25.6.25.1}).
We have the subcomplexes
$\nbigf^j_{\lambda}\gbigc_1(\nbigm)
=\lambda^j\gbigc_1(\nbigm)$.
We may naturally regard
$\gbigc_1(\nbigm)$
as a subcomplex of 
$\nbigmtilde\otimes\nbigl(t)\otimes
\Omegatilde^{\bullet}
_{\cnum_{\lambda}\times X\times\cnum_t/\cnum_{\lambda}}$.

\begin{lem}
\label{lem;24.6.25.3}
The inclusion
$\gbigc_1(\nbigm)\to
\nbigmtilde\otimes\nbigl(t)\otimes
\Omegatilde^{\bullet}
_{\cnum_{\lambda}\times X\times\cnum_t/\cnum_{\lambda}}$
is a quasi-isomorphism.
It induces quasi-isomorphisms
of subcomplexes
\[
 \nbigf_{\lambda}^k
 \gbigc_1(\nbigm)\lrarr
 \nbigf_{\lambda}^k\bigl(
\nbigmtilde\otimes\nbigl(t)\otimes
\Omegatilde^{\bullet}
 _{\cnum_{\lambda}\times X\times\cnum_t/\cnum_{\lambda}}
 \bigr)
\] 
and the quotient complexes:
\[
 \Gr^k_{\nbigf_{\lambda}}
 \gbigc_1(\nbigm)\lrarr
 \Gr^k_{\nbigf_{\lambda}}\bigl(
\nbigmtilde\otimes\nbigl(t)\otimes
\Omegatilde^{\bullet}
 _{\cnum_{\lambda}\times X\times\cnum_t/\cnum_{\lambda}}
 \bigr).
\]
\end{lem}
\pf
For any $a>-1$,
the induced morphisms
\[
 \Gr^V_{a}(\nbigmtilde)
 \to
 \lambda^{-1}
 \Gr^V_{a+1}(\nbigmtilde)\otimes dt
 \quad
 s\longmapsto \del_t(s)\otimes dt
\]
are isomorphisms.
We obtain that the quotient of
$\gbigc_1(\nbigm)\to
\nbigmtilde\otimes\nbigl(t)\otimes
\Omegatilde^{\bullet}
_{\cnum_{\lambda}\times X\times\cnum_t/\cnum_{\lambda}}$
is acyclic,
that is the first claim.
The second claim is equivalent to the first claim.
The third claim follows from the second claim.
\hfill\qed

\vspace{.1in}

By Lemma \ref{lem;24.6.25.3},
Theorem \ref{thm;23.3.29.2} is reduced to the following theorem.
\begin{thm}
\label{thm;23.3.29.3}
 The $E_1$-degeneration property holds for
the filtered complex
 $\nbigf_{\lambda}^{\bullet}
 \gbigc_1(\nbigm)$
with respect to the push-forward
by the projection
$\cnum_{\lambda}\times X\times\cnum_t
\to \cnum_{\lambda}$.
\end{thm}

\subsubsection{Refinement}

Let $\pi_{012}:\cnum_{\lambda}\times X\times\cnum_t\times\cnum_u
\to \cnum_{\lambda}\times X\times\cnum_t$
and
$\pi_{01}:\cnum_{\lambda}\times X\times\cnum_t\times\cnum_u
\to \cnum_{\lambda}\times X$
denote the projections.
We obtain the following complex
on $\cnum_{\lambda}\times X\times\cnum_t\times\cnum_u$:
\begin{equation}
\label{eq;23.4.1.1}
 \pi_{012}^{\ast}V_{-1}(\nbigmtilde)
 \lrarr
 \lambda^{-1}\pi_{012}^{\ast}V_0(\nbigmtilde)\otimes dt,
 \quad
 s\longmapsto (u\del_t+\lambda^{-1})s\,dt.
\end{equation}
The first term sits in the degree $0$.
We extend it to the following double complex:
\[
 \pi_{012}^{\ast}
 V_{-1}(\nbigmtilde)
 \otimes
 \pi_{01}^{\ast}
 \Omegatilde^{\bullet}_{\nbigx/\cnum}
 \lrarr
 \lambda^{-1}\pi_{012}^{\ast}
 V_0(\nbigmtilde)
 \otimes dt
 \otimes
 \pi_{01}^{\ast}
 \Omegatilde^{\bullet}_{\nbigx/\cnum}.
\]
Let $\gbigc(\nbigm)$ denote the associated complex
on $\cnum_{\lambda}\times X\times\cnum_t\times\cnum_u$.
We have the subcomplexes
$\nbigf_{\lambda}^j\gbigc(\nbigm)=\lambda^j\gbigc(\nbigm)$.

We consider the $\cnum^{\ast}$-action
on $\cnum_{\lambda}\times \cnum_u$
given by
$a(\lambda,u)=(a^{-1}\lambda,au)$.
It induces a $\cnum^{\ast}$-action
on $\cnum_{\lambda}\times X\times\cnum_t\times\cnum_u$.
Because $\nbigm$ is induced by a Hodge module,
$\pi_{012}^{\ast}V_{-1}(\nbigmtilde)$
and
$\lambda^{-1}\pi_{012}^{\ast}V_0(\nbigmtilde)$
are naturally $\cnum^{\ast}$-equivariant.
The differential
$u\del_t+\lambda^{-1}$ is $\cnum^{\ast}$-equivariant.

Let $\pi_{03}:\cnum_{\lambda}\times X\times\cnum_t\times\cnum_u
\lrarr \cnum_{\lambda}\times\cnum_u$
denote the projection.
The direct image sheaves
$R^{\ell}\pi_{03\ast}\gbigc(\nbigm)$
are $\cnum^{\ast}$-equivariant.

\begin{thm}
\label{thm;23.3.29.4}
The $E_1$-degeneration holds for
the filtered complex 
$\nbigf_{\lambda}^{\bullet}\gbigc(\nbigm)$
with respect to the push-forward by $\pi_{03}$,
and
$R^{\ell}\pi_{03\ast}
\Gr^j_{\nbigf_{\lambda}}\gbigc(\nbigm)$
are locally free
$\cnum^{\ast}$-equivariant $\nbigo_{\{0\}\times \cnum_u}$-modules
of finite rank.
\end{thm}

We obtain Theorem \ref{thm;23.3.29.3}
from Theorem \ref{thm;23.3.29.4}
by specializing along $u=1$.

\subsection{Formal neighbourhood along $u=0$}
\subsubsection{Basic strictness}

Let $\pi_0:\cnum_{\lambda}\times X\times\cnum_t\to \cnum_{\lambda}$
denote the projection.
Let $\iota_{u}:\cnum_{\lambda}\times X\times\cnum_t\times\{0\}\to
\cnum_{\lambda}\times X\times\cnum_t\times\cnum_u$
and
$\iota_{01}:\cnum_{\lambda}\times X\times\{0\}
 \lrarr
\cnum_{\lambda}\times X\times\cnum_t$
denote the inclusion maps.

\begin{lem}
\label{lem;23.3.29.10}
$R^{\ell}\pi_{0\ast}\bigl(
\iota_{u}^{\ast}\gbigc(\nbigm)
\bigr)$ are locally free
$\nbigo_{\cnum_{\lambda}}$-modules of finite rank
for any $\ell$.
\end{lem}
\pf
Because $\Gr^{V}_a(\nbigmtilde)$ $(-1<a\leq 0)$
underlie regular mixed twistor $\nbigd$-modules
on a K\"ahler manifold $X$
whose supports are compact,
the sheaves
$R^j\pi_{0\ast}\bigl(
\Gr^{V}_a(\nbigmtilde)
\otimes
\Omegatilde^{\bullet}_{\nbigx/\cnum}
\bigr)$
are locally free $\nbigo_{\cnum_{\lambda}}$-modules
of finite rank.
(See Corollary \ref{cor;25.6.29.1}.)

Note that
$\iota_u^{\ast}\gbigc(\nbigm)$
is quasi-isomorphic to
\[
       \lambda^{-1}\iota_{01\ast}\Bigl(
       \bigl(V_0(\nbigmtilde)/V_{-1}(\nbigmtilde)\bigr)
       \otimes
       dt\otimes
       \Omegatilde^{\bullet}_{\nbigx/\cnum}
       \Bigr)[-1].
\]       
For any $b<a$,
there exists the following canonical splitting of
$V\nbigr_{X\times\cnum_t|
\cnum_{\lambda}^{\ast}\times X\times\cnum_t}$-modules
\begin{equation}
\label{eq;23.3.29.10}
\bigl(V_a(\nbigmtilde)/V_{b}(\nbigmtilde)\bigr)
_{|\cnum_{\lambda}^{\ast}\times X\times\cnum_t}
=\bigoplus_{b<c\leq a}
\Gr^V_c(\nbigmtilde)_{|\cnum_{\lambda}^{\ast}\times X\times \cnum_t}.
\end{equation}
Here, the actions of $-\del_tt-c$
on 
$\Gr^V_c(\nbigmtilde)_{|\cnum_{\lambda}^{\ast}\times X\times \cnum_t}$
are nilpotent.
For any $-1<b<a\leq 0$,
we obtain the following exact sequence for any $\ell$:
\begin{multline}
 0\lrarr
 R^{\ell}\pi_{0\ast}\bigl(
  \Gr^V_b(\nbigmtilde)
  \otimes\Omegatilde^{\bullet}_{\nbigx/\cnum_{\lambda}}
  \bigr)
 _{|\cnum_{\lambda}^{\ast}}
\lrarr
 R^{\ell}\pi_{0\ast}\Bigl(
 \bigl(V_a(\nbigmtilde)/V_{<b}(\nbigmtilde)\bigr)
  \otimes\Omegatilde^{\bullet}_{\nbigx/\cnum_{\lambda}}
    \Bigr)
 _{|\cnum_{\lambda}^{\ast}}
 \\
  \lrarr
   R^{\ell}\pi_{0\ast}\Bigl(
   \bigl(V_a(\nbigmtilde)/V_{b}(\nbigmtilde)\bigr)
     \otimes\Omegatilde^{\bullet}_{\nbigx/\cnum_{\lambda}}
   \Bigr)
  _{|\cnum_{\lambda}^{\ast}}
  \lrarr 0.
\end{multline}
Because 
$R^{k}\pi_{0\ast}\bigl(
\Gr^V_b(\nbigmtilde)
\otimes\Omegatilde^{\bullet}_{\nbigx/\cnum_{\lambda}}
\bigr)$
are locally free $\nbigo_{\cnum_{\lambda}}$-modules
for any $k$,
we obtain the vanishing of the following morphisms:
\[
   R^{\ell}\pi_{0\ast}\Bigl(
   \bigl(V_a(\nbigmtilde)/V_{b}(\nbigmtilde)\bigr)
     \otimes\Omegatilde^{\bullet}_{\nbigx/\cnum_{\lambda}}
   \Bigr)
   \lrarr
    R^{\ell+1}\pi_{0\ast}\bigl(
  \Gr^V_b(\nbigmtilde)
  \otimes\Omegatilde^{\bullet}_{\nbigx/\cnum_{\lambda}}
  \bigr).
\]
We obtain the following exact sequences:
 \begin{multline}
 0\lrarr
 R^{\ell}\pi_{0\ast}\bigl(
  \Gr^V_b(\nbigmtilde)
  \otimes\Omegatilde^{\bullet}_{\nbigx/\cnum_{\lambda}}
  \bigr)
\lrarr
 R^{\ell}\pi_{0\ast}\Bigl(
 \bigl(V_a(\nbigmtilde)/V_{<b}(\nbigmtilde)\bigr)
  \otimes\Omegatilde^{\bullet}_{\nbigx/\cnum_{\lambda}}
    \Bigr)
 \\
  \lrarr
   R^{\ell}\pi_{0\ast}\Bigl(
   \bigl(V_a(\nbigmtilde)/V_{b}(\nbigmtilde)\bigr)
     \otimes\Omegatilde^{\bullet}_{\nbigx/\cnum_{\lambda}}
   \Bigr)
  \lrarr 0.
  \end{multline}  
Then, we obtain the claim of Lemma \ref{lem;23.3.29.10}
by an easy induction.
\hfill\qed

\subsubsection{The $E_1$-degeneration for another filtration
$\nbigf_u$}

We consider the subcomplexes
$\nbigf_u^j=u^j\gbigc(\nbigm)$.
There exist natural isomorphisms
$\Gr_{\nbigf_u}^j\gbigc(\nbigm)
\simeq
\iota_{u\ast}\iota_{u}^{\ast}\gbigc(\nbigm)$.
There exists the following exact sequence:
\begin{equation}
0\lrarr
      \Gr_{\nbigf_u}^j\gbigc(\nbigm)
      \lrarr
      \gbigc(\nbigm)/\nbigf_u^{j+1}\gbigc(\nbigm)
 \lrarr
      \gbigc(\nbigm)/\nbigf_u^{j}\gbigc(\nbigm)
\lrarr 0.
\end{equation}
\begin{prop}
\label{prop;23.3.29.30}
We obtain the following exact sequences:
\begin{equation}
 \label{eq;23.3.29.13}
0\lrarr
      R^{\ell}\pi_{03\ast}
      \Gr_{\nbigf_u}^j\gbigc(\nbigm)
 \lrarr
            R^{\ell}\pi_{03\ast}\bigl(
      \gbigc(\nbigm)/\nbigf_u^{j+1}\gbigc(\nbigm)
 \bigr)
 \\
      \lrarr
            R^{\ell}\pi_{03\ast}\bigl(
      \gbigc(\nbigm)/\nbigf_u^{j}\gbigc(\nbigm)
      \bigr)
\lrarr0.
 \end{equation}
\end{prop}
\pf
Let us consider the following exact sequence:
\[
      0\lrarr
 \nbigf_u^1\gbigc(\nbigm)\big/
 \nbigf_u^j\gbigc(\nbigm)
 \lrarr
 \gbigc(\nbigm)/\nbigf_u^{j}\gbigc(\nbigm)
 \\
      \lrarr
      \Gr_{\nbigf_u}^0\gbigc(\nbigm)
\lrarr0.
\]

\begin{lem}
\label{lem;23.3.29.12}
We obtain the following exact sequences
for any $\ell\in\seisuu_{\geq 0}$
and $j\in\seisuu_{\geq 1}$:
\begin{multline}
      0\lrarr
 R^{\ell}\pi_{03\ast}
 \bigl(
 \nbigf_u^1\gbigc(\nbigm)\big/
 \nbigf_u^j\gbigc(\nbigm)
 \bigr)
 _{|\cnum^{\ast}_{\lambda}\times\cnum_u}
 \lrarr
 R^{\ell}\pi_{03\ast}
 \bigl(
 \gbigc(\nbigm)/\nbigf_u^{j}\gbigc(\nbigm) \bigr)
 _{|\cnum^{\ast}_{\lambda}\times\cnum_u}
 \\
      \lrarr
 R^{\ell}\pi_{03\ast}\bigl(
 \Gr_{\nbigf_u}^0\gbigc(\nbigm)
 \bigr)
 _{|\cnum^{\ast}_{\lambda}\times\cnum_u}
\lrarr 0.
\end{multline}
\end{lem}
\pf
Recall that
$\Gr^0_{\nbigf_u}\gbigc(\nbigm)
=\iota_{u\ast}\iota_u^{\ast}\gbigc(\nbigm)$
is quasi-isomorphic to
\[
 \iota_{u\ast}\Bigl(
 \lambda^{-1}\iota_{01\ast}
 \Bigl(
 \pi_{012}^{\ast}
 \bigl(V_0(\nbigmtilde)/V_{-1}(\nbigmtilde)\bigr)
 \otimes dt\otimes
 \pi_{01}^{\ast}
 \Omegatilde^{\bullet}_{\nbigx/\cnum_{\lambda}}
 \Bigr)[-1]
 \Bigr).
\]
Let $\pi_1:\cnum_{\lambda}\times X\times\cnum_t\times\cnum_u\to X$
denote the projection.
Let $\nbigc_2^{\bullet}$
denote
the Dolbeault resolution,
i.e.,
the associated complex of the double complex
\[
  \iota_{u\ast}\Bigl(
 \lambda^{-1}\iota_{01\ast}
 \Bigl(
 \pi_{012}^{\ast}
 \bigl(V_0(\nbigmtilde)/V_{-1}(\nbigmtilde)\bigr)
 \otimes dt\otimes
 \pi_{01}^{\ast}
 \Omegatilde^{\bullet}_{\nbigx/\cnum_{\lambda}}
 \Bigr)[-1]
 \Bigr)
 \otimes
 \pi_1^{\ast}\Omega_X^{0,\bullet}.
\]
For any $N\geq 0$,
there exists the following subcomplex
$V_{-N-1}\gbigc(\nbigm)$ of $\gbigc(\nbigm)$:
\[
 \pi_{012}^{\ast}V_{-N-1}(\nbigmtilde)
 \otimes
 \pi_{01}^{\ast}\Omegatilde^{\bullet}_{\nbigx/\cnum}
 \lrarr
 \lambda^{-1}
 \pi_{012}^{\ast}V_{-N}(\nbigmtilde)
 \otimes dt\otimes
 \pi_{01}^{\ast}\Omegatilde^{\bullet}_{\nbigx/\cnum}.
\]

Let $U$ be an open subset of $\cnum_{\lambda}^{\ast}\times\cnum_u$.
Let $\alpha$ be a section of $\nbigc_2^{\bullet}$
on $U\times X\times\cnum_t$ such that $d\alpha=0$.
Let $N>j+10$.
By using the splitting (\ref{eq;23.3.29.10}),
we can construct a section $\alphatilde$
of
$\Bigl(
\lambda^{-1}\pi_{012}^{\ast}V_{0}(\nbigmtilde)
\otimes dt
\otimes
\pi_{01}^{\ast}\Omegatilde^{\bullet}_{\nbigx/\cnum}[-1]
\Bigr)
\otimes\pi_1^{\ast}\Omega^{0,\bullet}_X$
such that
(i) $\alphatilde$ induces $\alpha$,
(ii)
$\beta^{(0)}:=d\alphatilde$
is contained in
$\Bigl(
\lambda^{-1}
\pi_{012}^{\ast}
V_{-N}(\nbigmtilde)\otimes dt\otimes
\pi_{01}^{\ast}\Omegatilde^{\bullet}_{\nbigx/\cnum}
\Bigr) \otimes\pi_1^{\ast}\Omega_X^{0,\bullet}$
on $U\times X\times\cnum_t$.
Note that $d\beta^{(0)}=0$.
There exists a section
$\gamma^{(1)}$
of
$\Bigl(
\pi_{012}^{\ast}
V_{-N}(\nbigmtilde)\otimes
\pi_{01}^{\ast}\Omegatilde^{\bullet}_{\nbigx/\cnum}
\Bigr) \otimes\pi_1^{\ast}\Omega_X^{0,\bullet}$
on $U\times X\times\cnum_t$
such that
$\gamma^{(1)}\otimes \lambda^{-1}dt=\beta^{(0)}$.
We set
$\beta^{(1)}=\beta^{(0)}-d\gamma^{(1)}$
which is contained in
\[
 u\Bigl(
\lambda^{-1}
\pi_{012}^{\ast}
V_{-N+1}(\nbigmtilde)\otimes dt\otimes
\pi_{01}^{\ast}\Omegatilde^{\bullet}_{\nbigx/\cnum}
\Bigr) \otimes\pi_1^{\ast}\Omega_X^{0,\bullet}
\]
on $U\times X\times\cnum_t$.
We have $d\beta^{(1)}=0$.
Inductively,
for $m=1,\ldots,j+1$,
we can construct sections
$\beta^{(m)}$
of 
\[
u^m\Bigl(
\lambda^{-1}
\pi_{012}^{\ast}
V_{-N+m}(\nbigmtilde)\otimes dt\otimes
\pi_{01}^{\ast}\Omegatilde^{\bullet}_{\nbigx/\cnum}
\Bigr) \otimes \pi_1^{\ast}\Omega_X^{0,\bullet}
\]
and sections
$\gamma^{(m)}$
of
$u^{m-1}\Bigl(
\pi_{012}^{\ast}
V_{-N+m-1}(\nbigmtilde)\otimes
\pi_{01}^{\ast}\Omegatilde^{\bullet}_{\nbigx/\cnum}
\Bigr) \otimes\pi_1^{\ast}\Omega_X^{0,\bullet}$
such that
$\beta^{(m)}=\beta^{(m-1)}-d\gamma^{(m)}$
and
$\gamma^{(m)}\otimes \lambda dt=\beta^{(m-1)}$.
Then, 
$\hat{\alpha}=\alphatilde-\sum_{m=1}^{j+1} \gamma^{(m)}$
is a section of
$\gbigc(\nbigm)\otimes\pi_1^{\ast}\Omega_X^{0,\bullet}$
such that
(i) $\hat{\alpha}$ induces $\alpha$,
(ii) $d\hat{\alpha}$
is a section of
$\nbigf_u^{j+1}\gbigc(\nbigm)\otimes\pi_1^{\ast}\Omega_X^{0,\bullet}[1]$.
Then, we obtain
the claim of Lemma \ref{lem;23.3.29.12}.
\hfill\qed

\begin{lem}
We obtain the following exact sequences
for any $\ell,j\in\seisuu_{\geq 0}$.
\begin{multline}
\label{eq;24.6.25.2}
 0\lrarr
      R^{\ell}\pi_{03\ast}
 \Gr^{j}_{\nbigf_u}\gbigc(\nbigm)_{|\cnum^{\ast}_{\lambda}\times\cnum_u}
 \lrarr
            R^{\ell}\pi_{03\ast}\bigl(
 \gbigc(\nbigm)/\nbigf_u^{j+1}\gbigc(\nbigm) \bigr)
 _{|\cnum^{\ast}_{\lambda}\times\cnum_u}
 \\
      \lrarr
            R^{\ell}\pi_{03\ast}\bigl(
      \gbigc(\nbigm)/\nbigf_u^{j}\gbigc(\nbigm)
 \bigr)
 _{|\cnum^{\ast}_{\lambda}\times\cnum_u}
\lrarr0.
\end{multline}      
\end{lem}
\pf
It is enough to prove that
\begin{equation}
\label{eq;23.3.31.1}
             R^{\ell}\pi_{03\ast}\bigl(
 \gbigc(\nbigm)/\nbigf_u^{j+1}\gbigc(\nbigm) \bigr)
 _{|\cnum^{\ast}_{\lambda}\times\cnum_u}
      \lrarr
            R^{\ell}\pi_{03\ast}\bigl(
      \gbigc(\nbigm)/\nbigf_u^{j}\gbigc(\nbigm)
 \bigr)
 _{|\cnum^{\ast}_{\lambda}\times\cnum_u}
\end{equation}
are epimorphisms for any $\ell$ and $j$.
We use an induction on $j$.
Lemma \ref{lem;23.3.29.12} implies
the claim in the case $j=1$.
By the hypothesis of the induction,
we may assume that
\[
              R^{\ell}\pi_{03\ast}\bigl(
 \nbigf_u^1\gbigc(\nbigm)/\nbigf_u^{j+1}\gbigc(\nbigm) \bigr)
 _{|\cnum^{\ast}_{\lambda}\times\cnum_u}
      \lrarr
            R^{\ell}\pi_{03\ast}\bigl(
      \nbigf_u^1\gbigc(\nbigm)/\nbigf_u^{j}\gbigc(\nbigm)
 \bigr)
 _{|\cnum^{\ast}_{\lambda}\times\cnum_u}
\]
are epimorphisms.
Then, by Lemma \ref{lem;23.3.29.12},
we obtain that
(\ref{eq;23.3.31.1}) are also epimorphisms.
\hfill\qed

\vspace{.1in}
We obtain the exactness of (\ref{eq;23.3.29.13})
from the exactness of (\ref{eq;24.6.25.2})
and the strictness of
$R^{\ell}\pi_{03\ast}
\Gr_{\nbigf_u}^j\gbigc(\nbigm)$
for any $\ell$ and $j$
as in the proof of Lemma \ref{lem;23.3.29.10}.
Thus, we obtain Proposition \ref{prop;23.3.29.30}.
\hfill\qed

\vspace{.1in}
The supports of
$R^{\ell}\pi_{03\ast}\bigl(
\gbigc(\nbigm)/\nbigf_u^j\gbigc(\nbigm)
\bigr)$
are contained in
$\cnum_{\lambda}\times\{0\}\subset\cnum_{\lambda}\times\cnum_u$.
We may naturally regard them
as $\nbigo_{\cnum_{\lambda}}[u]/u^{j}$-modules.

\begin{cor}
\label{cor;23.4.1.2}
$R^{\ell}\pi_{03\ast}\bigl(
\gbigc(\nbigm)/\nbigf_u^j\gbigc(\nbigm)
\bigr)$
are locally free $\nbigo_{\cnum_{\lambda}}[u]/u^{j}$-modules.
\hfill\qed
\end{cor}

\subsection{Coherence of the specialization along $\lambda=0$}

Let $\iotatilde_{\lambda}:\{0\}\times X\times \cnum_t\times\cnum_u
\lrarr\cnum_{\lambda}\times X\times\cnum_t\times\cnum_u$
denote the inclusion.
Let $\pi_3:X\times\cnum_t\times\cnum_u\to\cnum_u$
denote the projection.

\begin{prop}
\label{prop;24.6.25.1}
$R^{\ell}\pi_{3\ast}(\iotatilde_{\lambda}^{\ast}\gbigc(\nbigm))$
are coherent $\nbigo_{\cnum_u}$-modules.
\end{prop}
\pf
Let $\pi_{12}:X\times\cnum_t\times\cnum_u\to X\times\cnum_t$
denote the projection.
Note that
$\iotatilde_{\lambda}^{\ast}\pi_{12}^{\ast}\nbigr_{X\times\cnum_t}$
is isomorphic to
the algebra of the symmetric product of
$\pi_{12}^{\ast}\Theta_{X\times\cnum_t}$,
where $\Theta_{X\times\cnum_t}$ denote the tangent sheaf of
$X\times\cnum_t$.
Because
$\iotatilde_{\lambda}^{\ast}\gbigc(\nbigm)$
is a complex of coherent
$\iotatilde_{\lambda}^{\ast}\nbigr_{X\times\cnum_t}$-modules,
it induces a complex of
coherent $\nbigo_{T^{\ast}(X\times\cnum_t)\times\cnum_u}$-modules
$\bigl(
\iotatilde_{\lambda}^{\ast}\gbigc(\nbigm)
\bigr)^{\sim}$.
Let $Z$ denote 
the cohomological support of
$\bigl(
\iotatilde_{\lambda}^{\ast}\gbigc(\nbigm)
\bigr)^{\sim}$.

\begin{lem}
Let $0_{X\times \cnum_t}:
X\times\cnum_t\to T^{\ast}(X\times\cnum_t)$
denote the $0$-section.
Then,
 $Z\subset
 \bigl(
 0_{X\times \cnum_t}\circ\iota_f\bigr)(\Cr(f))
 \times\cnum_u$.
\end{lem}
\pf
Let $\bigl(\iotatilde_{\lambda}^{\ast}\gbigc(\nbigm)\bigr)_u$
denote the pull back of
$\iotatilde_{\lambda}^{\ast}\gbigc(\nbigm)$
by the inclusion
$X\times \cnum_t\times\{u\}
\to X\times\cnum_t\times\cnum_u$.
Let $\bigl(
\iotatilde_{\lambda}^{\ast}\gbigc(\nbigm)
\bigr)_u^{\sim}$
denote the induced complex of
coherent $\nbigo_{T^{\ast}(X\times\cnum_t)}$-modules.
It equals the pull back of
$\bigl(
\iotatilde_{\lambda}^{\ast}\gbigc(\nbigm)
\bigr)^{\sim}$
by the inclusion
$T^{\ast}(X\times\cnum_t)\times\{u\}\to
T^{\ast}(X\times\cnum_t)\times\cnum_u$.
Let $Z_u$ denote the cohomological support of
$\bigl(
\iotatilde_{\lambda}^{\ast}\gbigc(\nbigm)
\bigr)_u^{\sim}$.
Because
$Z\cap\bigl(
T^{\ast}(X\times\cnum_t)\times\{u\}\bigr)
=Z_u$,
it is enough to prove that
$Z_u\subset (0_{X\times\cnum_t}\circ\iota_f)(\Cr(f))$.

Let $\iota_{0,\lambda}:\{0\}\times X\times\cnum_t\to
\cnum_{\lambda}\times X\times\cnum_t$ denote the inclusion.
By Lemma \ref{lem;24.6.25.3},
if $u\neq 0$,
$\bigl(
\iotatilde_{\lambda}^{\ast}\gbigc(\nbigm)
\bigr)_u$
is quasi-isomorphic to
$\iota_{0,\lambda}^{\ast}\Bigl(
\nbigmtilde\otimes\nbigl(u^{-1}t)\otimes
\Omegatilde^{\bullet}
_{\cnum_{\lambda}\times X\times\cnum_t/\cnum_{\lambda}}
\Bigr)$.
It is quasi-isomorphic to
$\iota_{f\ast}\iota_{\lambda}^{\ast}
\Bigl(
 \nbigm_{u^{-1}f}\otimes
 \Omegatilde^{\bullet}_{\nbigx/\cnum_{\lambda}}[-1]
\Bigr)$.
Hence, by Lemma \ref{lem;24.6.25.10},
the cohomological support
of $\bigl(
\iotatilde_{\lambda}^{\ast}\gbigc(\nbigm)
\bigr)^{\sim}_u$
is contained in
$(0_{X\times\cnum_t}\circ\iota_f)(\Cr(f))$.

Let us consider the case $u=0$.
There exists the following quasi-isomorphism
\[
\bigl(
\iotatilde_{\lambda}^{\ast}\gbigc(\nbigm)
\bigr)_0
\simeq
\iota_{\lambda^{\ast}}
\bigl(
V_0(\nbigmtilde)/V_{-1}(\nbigmtilde)
\bigr)\otimes dt
\otimes\Omega^{\bullet}_{X}[-1].
\]
Because $df(X)\cap \Ch(\Xi_{\DR}(\nbigm))\subset
0_X(\Cr(f))$,
the support of
the $\Sym\Theta_X$-module
$\iota_{\lambda\ast}(V_0(\nbigmtilde)/V_{-1}(\nbigmtilde))$
is contained in $\iota_f(\Cr(f))$
as in the case of Lemma \ref{lem;25.8.26.10}.
By a similar argument to the proof of Lemma \ref{lem;24.6.25.10},
we obtain that
the cohomological support of
$\bigl(
\iotatilde_{\lambda}^{\ast}\gbigc(\nbigm)
\bigr)_0$
is contained in
$(0_{X\times\cnum_t}\circ\iota_f)(\Cr(f))$.
\hfill\qed

\vspace{.1in}
Because $Z$ is proper over $\cnum_u$,
we obtain the claim of Proposition  \ref{prop;24.6.25.1}.
\hfill\qed

\vspace{.1in}
Note that
$\Gr^j_{\nbigf_{\lambda}}\gbigc(\nbigm)$
are isomorphic to
$\iota_{\lambda\ast}\iota_{\lambda}^{\ast}\gbigc(\nbigm)$
for any $j\in\seisuu_{\geq 0}$.
The supports of
$R^{\ell}\pi_{03\ast}(\gbigc(\nbigm)/\nbigf^j_{\lambda}\gbigc(\nbigm))$
are contained in $\{0\}\times\cnum_u$.
We may naturally regard
$R^{\ell}\pi_{03\ast}(\gbigc(\nbigm)/\nbigf^j_{\lambda}\gbigc(\nbigm))$
as $\nbigo_{\cnum_u}$-modules.

\begin{cor}
$R^{\ell}\pi_{03\ast}\bigl(
\gbigc(\nbigm)/\nbigf_{\lambda}^j\gbigc(\nbigm)
\bigr)$
are coherent $\nbigo_{\cnum_u}$-modules
for any $\ell$ and $j$. 
\hfill\qed
\end{cor}

\subsection{Proof of Theorem \ref{thm;23.3.29.4}}

There exist the following exact sequences:
\begin{multline}
  0\lrarr
 \lambda^j\bigl(\gbigc(\nbigm)/\nbigf_u^k\gbigc(\nbigm)\bigr)
  \Big/\lambda^{j+1}\bigl(\gbigc(\nbigm)/\nbigf_u^k\gbigc(\nbigm)\bigr)
  \\
\lrarr
  \bigl(\gbigc(\nbigm)/\nbigf_u^k\gbigc(\nbigm)\bigr)
  \Big/\lambda^{j+1}\bigl(\gbigc(\nbigm)/\nbigf_u^k\gbigc(\nbigm)\bigr)
  \\
\lrarr 
 \bigl(\gbigc(\nbigm)/\nbigf_u^k\gbigc(\nbigm)\bigr)
  \Big/\lambda^j\bigl(\gbigc(\nbigm)/\nbigf_u^k\gbigc(\nbigm)\bigr)
\lrarr 0.
\end{multline}
We obtain the following morphisms of the stalks at
$(0,0)\in\cnum_{\lambda}\times\cnum_u$:
\begin{multline}
\label{eq;23.3.29.41}
 R^{\ell}\pi_{03\ast}\Bigl(
  \bigl(\gbigc(\nbigm)/\nbigf_u^k\gbigc(\nbigm)\bigr)
  \Big/\lambda^j\bigl(\gbigc(\nbigm)/\nbigf_u^k\gbigc(\nbigm)\bigr)
 \Bigr)_{(0,0)}
 \lrarr \\
  R^{\ell+1}\pi_{03\ast}\Bigl(
  \lambda^j\bigl(\gbigc(\nbigm)/\nbigf_u^k\gbigc(\nbigm)\bigr)
  \Big/
   \lambda^{j+1}\bigl(\gbigc(\nbigm)/\nbigf_u^k\gbigc(\nbigm)\bigr)
  \Bigr)_{(0,0)}.
\end{multline}
\begin{lem}
\label{lem;23.4.1.3}
The morphisms {\rm(\ref{eq;23.3.29.41})} are $0$.
\end{lem}
\pf
Because
the multiplication of $\lambda^j$ on
$R^{\ell}\pi_{03\ast}\bigl(
\gbigc(\nbigm)/\nbigf_u^k\gbigc(\nbigm)
\bigr)$
is a monomorphism by Corollary \ref{cor;23.4.1.2},
the induced morphism 
\[
R^{\ell}\pi_{03\ast}\Bigl(
      \bigl(\gbigc(\nbigm)/\nbigf_u^k\gbigc(\nbigm)\bigr)
      \big/\lambda^{j}
      \bigl(\gbigc(\nbigm)/\nbigf_u^k\gbigc(\nbigm)\bigr)
      \Bigr)_{(0,0)}
      \lrarr
R^{\ell+1}\pi_{03\ast}\Bigl(
      \lambda^j\bigl(\gbigc(\nbigm)/\nbigf_u^k\gbigc(\nbigm)\bigr)
      \Bigr)_{(0,0)}
\]      
is $0$. Hence, the morphism (\ref{eq;23.3.29.41}) is $0$.
\hfill\qed

\vspace{.1in}
There exist the following exact sequences:
\[
  0 \lrarr
  \Gr_{\nbigf_{\lambda}}^j\gbigc(\nbigm)
  \lrarr
  \gbigc(\nbigm)/\nbigf^{j+1}_{\lambda}\gbigc(\nbigm)
  \lrarr
  \gbigc(\nbigm)/\nbigf^{j}_{\lambda}\gbigc(\nbigm)
  \lrarr 0.
\]
We obtain the following induced morphisms of the stalks:
\begin{equation}
\label{eq;23.3.29.42}
   R^{\ell}\pi_{03\ast}\bigl(
  \gbigc(\nbigm)/\nbigf^j_{\lambda}\gbigc(\nbigm)
  \bigr)_{(0,0)}
\lrarr
  R^{\ell+1}\pi_{03\ast}\bigl(
  \Gr_{\nbigf_{\lambda}}^j\gbigc(\nbigm)
  \bigr)_{(0,0)}.
\end{equation}
\begin{lem}
The morphisms {\rm(\ref{eq;23.3.29.42})} are $0$.
\end{lem}
\pf
There exist the following commutative diagrams
of stalks at $(\lambda,u)=(0,0)$ for any $k\geq 0$:
{\footnotesize
\[
 \begin{CD}
  R^{\ell}\pi_{03\ast}\bigl(
  \gbigc(\nbigm)/\nbigf^j_{\lambda}\gbigc(\nbigm)
  \bigr)_{(0,0)}
  @>{b^{\ell}}>>
  R^{\ell+1}\pi_{03\ast}\bigl(
  \Gr_{\nbigf_{\lambda}}^j\gbigc(\nbigm)
  \bigr)_{(0,0)}\\
  @VVV @VVV\\
  R^{\ell}\pi_{03\ast}\bigl(
  (\gbigc(\nbigm)/\nbigf_u^k\gbigc(\nbigm))
  \big/\lambda^j(\gbigc(\nbigm)/\nbigf_u^k\gbigc(\nbigm))
  \bigr)_{(0,0)}
  @>{0}>>
  R^{\ell+1}\pi_{03\ast}\bigl(
  \lambda^j(\gbigc(\nbigm)/\nbigf_u^k\gbigc(\nbigm))
  \big/
   \lambda^{j+1}(\gbigc(\nbigm)/\nbigf_u^k\gbigc(\nbigm))
  \bigr)_{(0,0)}.\\\\
 \end{CD}
\]}
The lower horizontal arrow is $0$
by Lemma \ref{lem;23.4.1.3}.
By the construction,
$\gbigc(\nbigm)$ is flat over
$\nbigo_{\cnum_{\lambda}\times\cnum_u}$,
and hence
$\lambda^p\gbigc(\nbigm)\cap
 u^q\gbigc(\nbigm)
=\lambda^pu^q\gbigc(\nbigm)$ for any $p,q\in\seisuu_{\geq 0}$.
The right vertical arrows are identified
as follows:{\footnotesize
\[
\begin{CD}
 R^{\ell+1}\pi_{03\ast}
 \iota_{\lambda\ast}\iota_{\lambda}^{\ast}
 \gbigc(\nbigm)_{(0,0)}
 @>{\simeq}>>
R^{\ell+1}\pi_{03\ast}\bigl(
  \Gr_{\nbigf_{\lambda}}^j\gbigc(\nbigm)
 \bigr)_{(0,0)}
 \\
@VVV @VVV \\  
  R^{\ell+1}\pi_{03\ast}
 \iota_{\lambda\ast}\iota_{\lambda}^{\ast}\bigl(
 \gbigc(\nbigm)/\nbigf_u^k\gbigc(\nbigm)
 \bigr)_{(0,0)}
 @>{\simeq}>>
 R^{\ell+1}\pi_{03\ast}\bigl(
  \lambda^j(\gbigc(\nbigm)/\nbigf_u^k\gbigc(\nbigm))
  \big/
   \lambda^{j+1}(\gbigc(\nbigm)/\nbigf_u^k\gbigc(\nbigm))
  \bigr)_{(0,0)}.
\end{CD}
\]      }
Under the identification,
the image of $b^{\ell}$
is contained in the image of
the morphism
\[
 R^{\ell+1}\pi_{03\ast}
 \iota_{\lambda\ast}\iota_{\lambda}^{\ast}
      \nbigf_u^k\gbigc(\nbigm)_{(0,0)}
      \lrarr
 R^{\ell+1}\pi_{03\ast}
 \iota_{\lambda\ast}\iota_{\lambda}^{\ast}
      \gbigc(\nbigm)_{(0,0)}
\]
for any $k\geq 0$.
It equals the image of the morphism
\[
 u^k:
 R^{\ell+1}\pi_{03\ast}
 \iota_{\lambda\ast}\iota_{\lambda}^{\ast}
      \gbigc(\nbigm)_{(0,0)}
      \lrarr
 R^{\ell+1}\pi_{03\ast}
 \iota_{\lambda\ast}\iota_{\lambda}^{\ast}
      \gbigc(\nbigm)_{(0,0)}
\]
for any $k\geq 0$.
Because
$R^{\ell+1}\pi_{03\ast}
 \iota_{\lambda\ast}\iota_{\lambda}^{\ast}
      \gbigc(\nbigm)$
is $\nbigo_{\cnum_{\lambda}}$-coherent by Proposition \ref{prop;24.6.25.1},
we obtain $b^{\ell}=0$.
\hfill\qed

\vspace{.1in}

For any $j\geq 0$,
there exists $\epsilon>0$ such that
the following morphism is an epimorphism
on $\{|u|<\epsilon\}$:
\[
     R^{\ell}\pi_{03\ast}\bigl(
     \gbigc(\nbigm)/\nbigf_{\lambda}^{j+1}\gbigc(\nbigm)
 \bigr)
      \lrarr \\
            R^{\ell}\pi_{03\ast}\bigl(
      \gbigc(\nbigm)/\nbigf_{\lambda}^j\gbigc(\nbigm)
      \bigr).
\]
By using the $\cnum^{\ast}$-equivariance,
we obtain that
it is an epimorphism on $\cnum_u$.
Thus, we obtain Theorem \ref{thm;23.3.29.4}.
\hfill\qed


\begin{thebibliography}{99}
\bibitem{Arinkin-Caldararu-Hablicsek}
	D. Arinkin,
	A. C\u{a}ld\u{a}raru,
	M. Hablicsek,
	{\em Derived intersections and the Hodge theorem},
	Algebr. Geom. {\bf 4} (2017), 394--423.

\bibitem{Esnault-Sabbah-Yu}
H. Esnault, 
C. Sabbah, 
J.-D. Yu,
(with an appendix by M. Saito),
{\em $E_1$-degeneration of the irregular Hodge filtration},
J. reine angew. Math. (2015), doi:10.1515/crelle-2014-0118.

\bibitem{Fan}
	 H. Fan,
	 {\em Schr\"{o}dinger equations, deformation theory
	 and $tt^{\ast}$-geometry},
	 arXiv:1107.1290

 \bibitem{Gelfand-Manin}
S. I. Gelfand, Y. I. Manin,
{\em Methods of homological algebra},
Springer-Verlag, Berlin, 1996, xviii+372 pp.
	
\bibitem{kashiwara_text}
M. Kashiwara,
{\em $D$-modules and microlocal calculus},
Translations of Mathematical Monographs, 217. 
Iwanami Series in Modern Mathematics,
American Mathematical Society, 
2003

\bibitem{Katzarkov-Kontsevich-Pantev-2017}
L. Katzarkov, M. Kontsevich, T. Pantev,
{\em Bogomolov-Tian-Todorov theorems for Landau-Ginzburg models},
J. Differential Geom. 105 (2017), no. 1, 55--117.

\bibitem{Li-Wen}
	S. Li,
	H. Wen,
	{\em On the $L^2$-Hodge theory of Landau-Ginzburg models},
	Adv. Math. {\bf 396} (2022), Paper No. 108165, 48 pp.
	
\bibitem{malgrange2}
B. Malgrange,
{\em Ideals of differentiable functions},
Tata Institute of Fundamental Research Studies in Mathematics
{\bf 3}, Tata Institute of Fundamental Research, Bombay,
Oxford University Press, London, 1967.	 
	
\bibitem{mochi2}
T. Mochizuki,
{\em Asymptotic behaviour of tame harmonic bundles
and an application to pure twistor $\nbigd$-modules I, II},
Mem. AMS. {\bf 185}, (2007).
	
\bibitem{Mochizuki-wild}
T. Mochizuki,
{\em Wild harmonic bundles and 
 wild pure twistor $\nbigd$-modules},
Ast\'{e}risque {\bf 340}, (2011)

\bibitem{Mochizuki-MTM}
T. Mochizuki,
{\em Mixed twistor $\nbigd$-modules},
Springer, 2015.

\bibitem{Mochizuki-Kontsevich-complexes}
	T. Mochizuki,
	{\em A twistor approach to the Kontsevich complexes},
	Manuscripta Math. {\bf 157} (2018), 193--231.

 \bibitem{Mochizuki-L2}
	 T. Mochizuki,
	 {\em $L^2$-complexes and twistor complexes of
	 tame harmonic bundles},
	 arXiv:2204.10443	 
 \bibitem{Ogus-Vologodsky}
	 A. Ogus, 
	 V. Vologodsky,
	 {\em Nonabelian Hodge theory in characteristic $p$.}
	 Publ. Math. Inst. Hautes \'{E}tudes Sci. {\bf 106} (2007), 1--138.
	 
 \bibitem{Sabbah-twisted-de-Rham}
	 C. Sabbah,
	 {\em On a twisted de Rham complex},
	 Tohoku Math. J. (2) {\bf 51} (1999), 125--140.
	 
\bibitem{Sabbah-pure-twistor}
	 C. Sabbah, {\em Polarizable twistor $\nbigd$-modules},
	 Ast\'{e}risque {\bf 300} (2005).

 \bibitem{Sabbah-tame-polynomial}
	 C. Sabbah,
	 {\em Hypergeometric periods for a tame polynomial},
	 Port. Math. (N.S.) {\bf 63} (2006), no. 2, 173--226.
	
\bibitem{Sabbah-wild-twistor}
	C. Sabbah,
	{\em Wild twistor $\nbigd$-modules},
	in {\em Algebraic analysis and around},
	293--353, Adv. Stud. Pure Math., {\bf 54},
	Math. Soc. Japan, Tokyo, 2009. 

 \bibitem{Sabbah-Brieskorn}
	 C. Sabbah,
	 {\em Some properties and applications of Brieskorn lattices,}
	 J. Singul. {\bf 18} (2018), 238--247.
 
 \bibitem{Saito-Takahashi}
	 K. Saito, A. Takahashi,
	 {\em From primitive forms to Frobenius manifolds},
	 Proc. Sympos. Pure Math., {\bf 78}
	 American Mathematical Society, Providence, RI, 2008, 31--48.
	
\bibitem{saito1}
M. Saito,
{\em Modules de Hodge polarisables},
Publ. RIMS., {\bf 24},
(1988), 849--995.

\bibitem{saito2}
M. Saito,
{\em Mixed Hodge modules},
Publ. RIMS., {\bf 26}, (1990),
221--333.

\bibitem{s3}
C. Simpson,
{\it Mixed twistor structures},
math.AG/9705006.
	
\end{thebibliography}
\end{document}